\documentstyle[amscd,amssymb,xypic,verbatim,11pt]{amsart}
\xyoption{all}

\newcommand{\nc}{\newcommand}
\newcommand{\rnc}{\renewcommand}

\nc{\exto}[1]{\stackrel{#1}{\longrightarrow}}
\nc{\dlim}{{\mathop{\lim\limits_{\longrightarrow}}}}
\nc{\lan}{\big\langle}
\nc{\ran}{\big\rangle}

\nc{\kk}{{\mathsf{k}}}
\nc{\ix}{{\mathsf{i}}}
\nc{\jx}{{\mathsf{j}}}

\nc{\C}{{\mathbb{C}}}
\nc{\HH}{{\mathbb{H}}}
\nc{\PP}{{\mathbb{P}}}
\nc{\QQ}{{\mathbb{Q}}}
\nc{\ZZ}{{\mathbb{Z}}}

\nc{\CA}{{\mathcal{A}}}
\nc{\CB}{{\mathcal{B}}}
\nc{\CC}{{\mathcal{C}}}
\nc{\D}{{\mathcal{D}}}
\nc{\CE}{{\mathcal{E}}}
\nc{\CF}{{\mathcal{F}}}
\nc{\CG}{{\mathcal{G}}}
\nc{\CH}{{\mathcal{H}}}
\nc{\CJ}{{\mathcal{J}}}
\nc{\CK}{{\mathcal{K}}}
\nc{\CL}{{\mathcal{L}}}
\nc{\CM}{{\mathcal{M}}}
\nc{\CN}{{\mathcal{N}}}
\nc{\CO}{{\mathcal{O}}}
\nc{\CQ}{{\mathcal{Q}}}
\nc{\CR}{{\mathcal{R}}}
\nc{\CS}{{\mathcal{S}}}
\nc{\CT}{{\mathcal{T}}}
\nc{\CU}{{\mathcal{U}}}
\nc{\CV}{{\mathcal{V}}}
\nc{\CW}{{\mathcal{W}}}
\nc{\CX}{{\mathcal{X}}}
\nc{\CY}{{\mathcal{Y}}}
\nc{\CZ}{{\mathcal{Z}}}
\nc{\CMo}{{\mathcal{M}^\circ}}
\nc{\Co}{{{C}^\circ}}

\nc{\BY}{{\overline{Y}}}
\nc{\BYD}{{\overline{Y}{}^{|D|}}}
\nc{\OZ}{{\overline{Z}}}
\nc{\bg}{{\bar{g}}}

\nc{\bq}{{\mathbf{q}}}
\nc{\BD}{{\mathbf{D}}}
\nc{\BG}{{\mathbf{G}}}
\nc{\BL}{{\mathbf{L}}}
\nc{\BM}{{\mathbf{M}}}
\nc{\BP}{{\mathbf{P}}}
\nc{\BZ}{{\mathbf{Z}}}
\nc{\BPr}{{\mathsf{P}}}
\nc{\BR}{{\mathbf{R}}}
\nc{\BRO}[1]{{{\mathbf{R}}^{\circ}_{#1}}}
\nc{\BRD}[1]{{{\mathbf{R}}^{|D|}_{#1}}}
\nc{\BRP}[1]{{{\mathbf{R}}^{1}_{#1}}}
\nc{\BRTP}[1]{{{\mathbf{\tilde{R}}}{}^{1}_{#1}}}
\nc{\BS}{{\mathbf{S}}}
\nc{\BMS}{{{\mathbf{M}}^{{s}}}}
\nc{\BMSS}{{{\mathbf{M}}^{{ss}}}}
\nc{\BMZ}{{\mathbf{M}^{\circ}}}
\nc{\BCL}{{\mathbf{L}}}

\nc{\PCC}{{{}^\perp\CC}}

\nc{\Cl}{{\mathsf{Cliff}}}
\nc{\Clev}{{\mathop{\mathsf{Cliff}}^{\circ}}}

\nc{\FA}{{\mathfrak{A}}}
\nc{\FB}{{\mathfrak{B}}}
\nc{\FI}{{\mathfrak{I}}}
\nc{\FZ}{{\mathfrak{Z}}}

\nc{\TFA}{{\tilde{\mathfrak{A}}}}
\nc{\TFB}{{\tilde{\mathfrak{B}}}}

\nc{\fa}{{\mathfrak{a}}}
\nc{\fg}{{\mathfrak{g}}}
\nc{\fp}{{\mathfrak{p}}}
\nc{\FD}{{\mathfrak{D}}}
\nc{\FE}{{\mathfrak{E}}}
\nc{\FL}{{\mathfrak{L}}}
\nc{\FM}{{\mathfrak{M}}}
\nc{\FR}{{\mathfrak{R}}}
\nc{\FS}{{\mathsf{S}}}

\nc{\sfc}{{\mathsf{c}}}
\nc{\sfch}{{\mathsf{ch}}}
\nc{\sfh}{{\mathsf{h}}}

\nc{\SK}{{\mathsf{K}}}
\nc{\SO}{{\mathsf{O}}}
\nc{\SQ}{{\mathsf{Q}}}
\nc{\SPV}{{\mathsf{S}^+\mathsf{V}}}
\nc{\SMV}{{\mathsf{S}^-\mathsf{V}}}
\nc{\SPMV}{{\mathsf{S}^\pm\mathsf{V}}}
\nc{\SX}{{S_X}}
\nc{\SY}{{S_Y}}
\nc{\phipsi}{{q}}
\nc{\eps}{\varepsilon}

\nc{\pim}{{\pi_-}}
\nc{\pip}{{\pi_+}}

\nc{\BE}{{\overline{\CE}}}
\nc{\TE}{{\tilde{\CE}}}
\nc{\TQ}{{\tilde{Q}}}
\nc{\TCF}{{\tilde{\CF}}}
\nc{\TCG}{{\tilde{\CG}}}
\nc{\TCL}{{\tilde{\CL}}}
\nc{\TF}{{\tilde{F}}}
\nc{\TW}{{\tilde{W}}}
\nc{\TCB}{{\widetilde{\CB}}}
\nc{\TCC}{{\tilde{\CC}}}
\nc{\TCX}{{\tilde{\CX}}}
\nc{\TCY}{{\tilde{\CY}}}
\nc{\TPhi}{{\tilde{\Phi}}}
\nc{\OPhi}{{\bar{\Phi}}}
\nc{\txi}{{\tilde{\xi}}}
\nc{\tp}{{\tilde{p}}}
\nc{\tq}{{\tilde{q}}}
\nc{\tzeta}{{\tilde{\zeta}}}
\nc{\tpi}{{\tilde{\pi}}}

\nc{\HCB}{{\widehat{\CB}}}
\nc{\HE}{{\widehat{\CE}}}
\nc{\HS}{{\widehat{S}}}
\nc{\HX}{{\hat{X}}}
\nc{\HY}{{\hat{Y}}}
\nc{\HZ}{{\hat{Z}}}
\nc{\hxi}{{\hat{\xi}}}

\nc{\UH}{{\mathcal{H}}}

\nc{\TD}{{\widetilde{\D}}}
\nc{\TM}{{\widetilde{M}}}
\nc{\TCM}{{\widetilde{\CM}}}
\nc{\TS}{{\widetilde{S}}}
\nc{\TT}{{\widetilde{T}}}
\nc{\TU}{{\widetilde{U}}}
\nc{\TX}{{\widetilde{X}}}
\nc{\TY}{{\widetilde{Y}}}
\nc{\TZ}{{\widetilde{Z}}}
\nc{\TYO}{{{\widetilde{Y}}^\circ}}
\nc{\barf}{{\bar{f}}}
\nc{\te}{{\tilde{e}}{}}
\nc{\tf}{{\tilde{f}}}
\nc{\tg}{{\tilde{g}}}
\nc{\ti}{{\tilde{\imath}}}
\nc{\tj}{{\tilde{\jmath}}}
\nc{\ty}{{\tilde{y}}}
\nc{\tphi}{{\tilde{\phi}}}
\nc{\hf}{{\mathsf{hf}}}
\nc{\hub}{{\mathsf{hub}}}

\nc{\urho}{{\underline{\rho}}}

\nc{\LRA}{\Leftrightarrow}
\nc{\RA}{\Rightarrow}
\nc{\lotimes}{\mathbin{\mathop{\otimes}\limits^{\mathbb{L}}}}
\nc{\CEnd}{\mathop{\mathcal{E}\mathit{nd}}\nolimits}
\nc{\CExt}{\mathop{\mathcal{E}\mathit{xt}}\nolimits}
\nc{\CHom}{\mathop{\mathcal{H}\mathit{om}}\nolimits}
\nc{\RH}{\mathop{{\mathsf{R}}\Gamma}\nolimits}
\nc{\RGamma}{\mathop{{\mathsf{R}}\Gamma}\nolimits}
\nc{\RHom}{\mathop{\mathsf{RHom}}\nolimits}
\nc{\RCHom}{\mathop{\mathsf{R}\mathcal{H}\mathit{om}}\nolimits}
\nc{\RG}{\mathop{\mathsf{R\Gamma}}\nolimits}
\nc{\Hom}{\mathop{\mathsf{Hom}}\nolimits}
\nc{\Ext}{\mathop{\mathsf{Ext}}\nolimits}
\nc{\End}{\mathop{\mathsf{End}}\nolimits}
\nc{\Tor}{\mathop{\mathsf{Tor}}\nolimits}
\nc{\Tordim}{\mathop{\mathsf{Tor}\text{\rm-}\mathsf{dim}}\nolimits}
\nc{\Hilb}{\mathop{\mathsf{Hilb}}\nolimits}
\nc{\Spec}{\mathop{\mathsf{Spec}}\nolimits}
\nc{\Proj}{\mathop{\mathsf{Proj}}\nolimits}
\nc{\Pic}{\mathop{\mathsf{Pic}}\nolimits}

\nc{\Tw}{\mathop{\mathsf{Tw}}\nolimits}
\nc{\Cone}{\mathop{\mathsf{Cone}}\nolimits}
\nc{\Ker}{\mathop{\mathsf{Ker}}\nolimits}
\nc{\Coker}{\mathop{\mathsf{Coker}}\nolimits}
\nc{\codim}{\mathop{\mathsf{codim}}\nolimits}
\nc{\discr}{{\mathsf{discr}}}
\nc{\sing}{{\mathsf{sing}}}
\nc{\supp}{\mathop{\mathsf{supp}}}
\nc{\perf}{{\mathsf{perf}}}
\nc{\rank}{\mathop{\mathsf{rank}}}
\nc{\Pf}{{\mathsf{Pf}}}
\nc{\Gr}{{\mathsf{Gr}}}
\nc{\OGr}{{\mathsf{OGr}}}
\nc{\Flag}{{\mathsf{Fl}}}
\nc{\Kosz}{{\mathsf{Kosz}}}
\nc{\LGr}{{\mathsf{LGr}}}
\nc{\GTGr}{{\mathsf{G_2Gr}}}
\nc{\GTF}{{\mathsf{G_2F}}}
\nc{\OF}{{\mathsf{OF}}}
\nc{\Fl}{{\mathsf{Fl}}}
\nc{\Bl}{{\mathsf{Bl}}}
\nc{\GL}{{\mathsf{GL}}}
\nc{\PGL}{{\mathsf{PGL}}}
\nc{\SL}{{\mathsf{SL}}}
\nc{\SP}{{\mathsf{Sp}}}
\nc{\Spin}{{\mathsf{Spin}}}
\nc{\Tot}{{\mathsf{Tot}}}
\nc{\ev}{{\mathsf{ev}}}
\nc{\od}{{\mathsf{odd}}}
\nc{\coev}{{\mathsf{coev}}}
\nc{\id}{{\mathsf{id}}}
\nc{\rk}{{\mathsf{r}}}
\nc{\opp}{{\mathsf{opp}}}
\nc{\PS}{{{\PP^3}}}
\nc{\Qu}{{{Q^3}}}
\nc{\tdim}{\mathop{\Tor\dim}}
\nc{\ecart}{{\fbox{$\scriptstyle\mathsf{EC}$}}}
\nc{\ad}{{\mathop{\mathsf ad}}}
\nc{\gr}{{\mathop{\mathsf gr}}}
\nc{\qgr}{{\mathop{\mathsf qgr}}}
\nc{\tor}{{\mathop{\mathsf tor}}}
\rnc{\mod}{{\mathop{\mathsf mod}}}
\nc{\Mod}{{\mathop{\mathsf Mod}}}
\nc{\Coh}{{\mathop{\mathsf Coh}}}
\nc{\Ab}{{\mathop{\mathcal{A}\mathit{b}}}}
\nc{\QCoh}{{\mathop{\mathsf QCoh}}}

\nc{\AAV}{{\mathcal{AAV}}}

\nc{\Rep}{{\mathsf{Rep}}}

\nc{\Cubics}{{{\mathcal{S}}_3}}
\nc{\VFT}{{{\mathcal{S}}_{14}}}
\nc{\VFTE}{{{\mathcal{N}}_{\mathrm{reg,sm}}}}
\nc{\MX}{{\CM_X}}
\nc{\MY}{{\CM_Y}}
\nc{\MYE}{{\CM_{Y,\CE}}}
\nc{\Yd}{{Y_d}}
\nc{\Yfive}{{Y_5}}
\nc{\Xg}{{X_{2g-2}}}
\nc{\Xtt}{{X_{22}}}
\nc{\Xst}{{X_{16}}}
\nc{\Xtw}{{X_{12}}}
\nc{\Xe}{{X_{8}}}
\nc{\Xf}{{X_{4}}}

\nc{\git}{{/\!\!/\!{}_\chi}}

\theoremstyle{plain}

\newtheorem{theo}{Theorem}[]
\newtheorem{theorem}{Theorem}[section]
\newtheorem{conjecture}[theorem]{Conjecture}
\newtheorem{lemma}[theorem]{Lemma}
\newtheorem{proposition}[theorem]{Proposition}
\newtheorem{corollary}[theorem]{Corollary}

\theoremstyle{definition}

\newtheorem{definition}[theorem]{Definition}

\newtheorem{example}[theorem]{Example}

\theoremstyle{remark}

\newtheorem{remark}[theorem]{Remark}

\newenvironment{proof}{\noindent{\sf Proof:}}{\qed\medskip}

\title[Lefschetz decompositions and categorical resolutions of singularities]%
{Lefschetz decompositions\\and\\categorical resolutions of singularities}
\author{Alexander Kuznetsov}
\address{
Algebra Section, Steklov Mathematical Institute,
8 Gubkin str., Moscow 119991 Russia
}
\email{akuznet@@mi.ras.ru}
\date{}
\thanks{I was partially supported by RFFI grants 05-01-01034 and 02-01-01041,
Russian Presidential grant for young scientists No. MK-6122.2006.1,
CRDF Award No. RUM1-2661-MO-05 and gratefully acknowledge the support
of the Pierre Deligne fund based on his 2004 Balzan prize in mathematics.}

\begin{document}

\begin{abstract}
Let $Y$ be a singular algebraic variety and let $\TY$ be a resolution of singularities of $Y$.
Assume that the exceptional locus of $\TY$ over $Y$ is an irreducible divisor $\TZ$ in $\TY$.
For every Lefschetz decomposition of $\TZ$ we construct a triangulated subcategory
$\TD \subset \D^b(\TY)$ which gives a desingularization of $\D^b(Y)$.
If the Lefschetz decomposition is generated by a vector bundle tilting over $Y$
then $\TD$ is a noncommutative resolution, and if the Lefschetz decomposition is rectangular,
then $\TD$ is a crepant resolution.
\end{abstract}

\maketitle

\section{Introduction}\label{intro}

The Minimal Model Program (MMP for short) is one of the most important branches
of birational geometry.
In the classical approach, there is an unpleasant but inevitable step,
enlargement of the category of smooth varieties by varieties with terminal
singularities.
From a modern point of view~\cite{BO2}, MMP can be considered as a minimization
of the derived category of coherent sheaves on an algebraic variety in a given
birational class.
In a contrast with the classical approach, the modern point of view
allows one to stay in the world of ``smooth'' triangulated categories
while searching for a minimal model.
This is achieved through the notion of a minimal categorical
resolution of singularities.

Let $Y$ be a singular algebraic variety. If $\pi:\TY \to Y$ is a resolution
of singularities then the derived categories of coherent sheaves on $\TY$ and $Y$
are related by the derived pushforward and derived pullback functors:
$$
\pi_*:\D^b(\TY) \to \D^b(Y),
\qquad\text{and}\qquad
\pi^*:\D^\perf(Y) \to \D^b(\TY).
$$
Here $\D^b$ stands for the bounded derived category of coherent sheaves,
and $\D^\perf$ stands for the category of perfect complexes.
The functors $\pi^*$ and $\pi_*$ are mutually adjoint
($\pi^*$ is left adjoint to $\pi_*$).
Moreover, if the singularities of $Y$ are rational,
then the composition $\pi_*\circ\pi^*$ is isomorphic to the identity functor.
Note also, that if $\pi$ is a \emph{crepant} resolution (i.e.\ the relative canonical
class is trivial) then $\pi^*$ is isomorphic to the right adjoint functor $\pi^!$ of $\pi_*$.

We suggest to take this structure for a definition of a categorical resolution of singularities.
A \emph{categorical resolution} of $\D^b(Y)$ should be a ``smooth'' triangulated category $\TD$
and a pair of functors
$$
\pi_*:\TD \to \D^b(Y),
\qquad\text{and}\qquad
\pi^*:\D^\perf(Y) \to \TD,
$$
such that $\pi^*$ is left adjoint to $\pi_*$
and the natural morphism of functors $\id_{\D^\perf} \to \pi_*\pi^*$
is an isomorphism.
A categorical resolution is \emph{crepant}, if $\pi^*$ is simultaneously
a right adjoint of $\pi_*$.

Unfortunately, the last condition in the definition of a categorical resolution
(which should be thought of as a way to express ``surjectivity'' of the functor $\pi_*$)
restricts us to consider only rational singularities. To include all singularities
in consideration one should weaken somehow this condition.

Certainly, the above definition relies on a notion of smoothness for triangulated categories,
for which as yet there is no generally adopted definition as well. In fact, there are
several approaches. One of them, considering saturatedness and $\Ext$-boundedness
as a definition of smoothness, is very convenient and adequate, but only
in projective (or at least proper) case.
Other approaches \cite{TV,KS} suggest to consider only those triangulated categories
which are equivalent to the derived categories of $A_\infty$ or DG-algebras
and define smoothness of a category in terms of homological properties
of the corresponding algebra.
However, we use another definition.
We consider as smooth only those triangulated categories which are equivalent to
admissible subcategories of bounded derived categories of smooth varieties.
On one hand, these categories should be smooth in any definition of smoothness,
and on the other hand, they are sufficient for our purposes.


A categorical resolution of singularities is called \emph{minimal}~\cite{BO2}
if it can be embedded as a full subcategory into any other categorical
resolution of the same singularity. This definition certainly is very
ineffective. Moreover, up to now there is no method to prove minimality
of a resolution. However, there is some evidence that at least \emph{crepant}
resolutions of singularities should be minimal.

First examples of categorical resolutions of singularities
were investigated in the context of the McKey correspondence.
If $V = \Spec R$ is a smooth affine variety and $\Gamma$ a finite
group acting on $V$ generically free, then the quotient
$V/\Gamma = \Spec R^\Gamma$ is singular and
the derived category of modules over the wreath product algebra
$R\#\Gamma$ provides a categorical resolution of $\D^b(V/\Gamma)$,
see~\cite{BKR,BeKa}.

First examples of categorical crepant resolutions of non-quotient singularities
were discovered by Van den Bergh \cite{V1,V2} under the name of \emph{noncommutative}
crepant resolutions. The resolutions were realized as derived categories of a sheaf
of noncommutative (hence the name) algebras with good homological properties
on the resolved variety. Such resolutions were constructed for threefold
terminal singularities which admit commutative crepant resolutions
and for cones over del Pezzo surfaces.

Another work in this area was done by Kaledin \cite{Ka}. He constructed
noncommutative crepant resolutions of symplectic singularities admitting
commutative crepant resolutions in all dimensions.

The main goal of the present paper is to construct a categorical resolution
of singularities in the following situation.
Let $\pi:\TY \to Y$ be a resolution of rational singularities.
Assume that the exceptional locus of $\pi$ is an irreducible divisor $\TZ \subset \TY$.
We show that in this situation categorical resolutions of $\D^b(Y)$
which can be embedded in $\D^b(\TY)$ are related to Lefschetz
decompositions of $\D^b(\TZ)$.

A \emph{Lefschetz decomposition} of the derived category $\D^b(X)$ of an algebraic variety $X$
with respect to a line bundle $L$ on $X$ is a chain of triangulated subcategories
$0 \subset \CB_{m-1} \subset \CB_{m-2} \subset \dots \subset \CB_1 \subset \CB_0 \subset \D^b(X)$
such that
$\D^b(X) = \lan \CB_{m-1}\otimes L^{1-m},\CB_{m-2}\otimes L^{2-m},\dots,\CB_1\otimes L^{-1},\CB_0 \ran$
is a semiorthogonal decomposition.
Lefschetz decompositions were introduced in~\cite{K2} in a context of Homological Projective Duality.
However, this notion is very important and useful in more general situation.
Relation to categorical resolutions is a confirmation of this.

Now we can formulate one of the main results of the paper.
Let $L = \CN^*_{\TZ/\TY}$ be the conormal bundle, and let
$$
D^b(\TZ) = \lan \CB_{m-1}\otimes L^{1-m},\CB_{m-2}\otimes L^{2-m},\dots,\CB_1\otimes L^{-1},\CB_0 \ran
$$
be a Lefschetz decomposition such that $p^*(\D^\perf(Z)) \subset \CB_0$ and $\CB_0$
is stable with respect to tensoring by $p^*(\D^\perf(Z))$, where
$Z = \pi(\TZ) \subset Y$ and $p = \pi_{|\TZ}:\TZ \to Z$.
Let $i:\TZ \to \TY$ be the embedding.
Denote by $\TD$ the full subcategory of $\D^b(\TY)$ consisting of objects
$F \in \D^b(\TY)$ such that $i^*F \in \CB_0 \subset \D^b(\TZ)$.
Then the functor $\pi^*:\D^\perf(Y) \to \D^b(\TY)$ factors through $\D^\perf(Y) \to \TD$
and the restriction of the functor $\pi_*$ to $\TD$ is its right adjoint.

\begin{theo}\label{theo1}
Triangulated category $\TD$ with functors $\pi_*$ and $\pi^*$ is
a categorical resolution of $\D^b(Y)$. Moreover, we have a semiorthogonal
decomposition
$$
\D^b(\TY) =
\lan i_*(\CB_{m-1}\otimes L^{1-m}),i_*(\CB_{m-2}\otimes L^{2-m}),\dots,i_*(\CB_1\otimes L^{-1}),\TD \ran.
$$
Finally, if $Y$ is Gorenstein, $\CB_{m-1} = \CB_{m-2} = \dots = \CB_1 = \CB_0$ and
$K_\TY = \pi^*K_{Y} + (m-1)\TZ$ then $\TD$ is a crepant categorical resolution of $\D^b(Y)$.
\end{theo}

Certainly, the most interesting are minimal categorical resolutions.
We introduce a notion of minimality for Lefschetz decompositions
and conjecture that the categorical resolution corresponding to a minimal
Lefschetz decomposition is minimal.

Another question addressed in the paper is which of categorical resolutions
constructed in Theorem~\ref{theo1} are noncommutative in the sense of Van den Bergh?
The answer is the following. Recall that a vector bundle $E$ on $\TY$ is called
tilting over $Y$ if the pushforward $\pi_*\CEnd E$ is a pure sheaf
(i.e.\ $R^{>0}\pi_*\CEnd E = 0$).

\begin{theo}\label{theo2}
Assume that there exists a vector bundle $E$ on $\TY$ such that
the category $\CB_0 \subset \D^b(\TZ)$ is generated by $i^*E$ and $E$ is tilting over $Y$.
Assume also that $\CB_0$ is admissible and $\CJ_\TZ = \pi^{-1}\CJ_Z\cdot\CO_\TY$.
Then the sheaf of algebras $\CA = \pi_*\CEnd E$ has finite homological dimension
and the category $\TD \cong \D^b(Y,\CA)$ is a noncommutative resolution of $\D^b(Y)$.
Moreover, if $Y$ is Gorenstein, $\CB_{m-1} = \CB_{m-2} = \dots = \CB_1 = \CB_0$ and
$K_\TY = \pi^*K_{Y} + (m-1)\TZ$ then $\D^b(Y,\CA)$ is a noncommutative crepant resolution of $\D^b(Y)$.
\end{theo}

We also give several examples of categorical and noncommutative resolutions of singularities
obtained by Theorems~\ref{theo1} and~\ref{theo2}, namely cones over smooth projective
varieties and Pfaffian varieties.

The paper is organized as follows.
In section~\ref{prelim} we recall the necessary background.
In section~\ref{defres} we discuss a definition of categorical resolutions of singularities.
In section~\ref{ld_cr} we prove Theorem~1 and formulate conjectures relating minimality
of a Lefschetz decomposition to minimality of the corresponding categorical resolution.
In section~\ref{ncr} we prove Theorem~2.
In section~\ref{funct} we investigate functoriality properties of resolutions constructed
in sections~\ref{ld_cr} and~\ref{ncr}.
In section~\ref{CONES} we discuss categorical resolution of singularities
of cones over smooth projective varieties.
Finally, in section~\ref{PF} we construct a noncommutative resolution of singularities
for some Pfaffian varieties.

\bigskip

{\bf Acknowledgements:}
I would like to thank A.~Bondal, D.~Kaledin, D.~Orlov and T.~Pantev for
very helpful discussions.

\section{Preliminaries}\label{prelim}

\subsection{Notation}

All algebraic varieties are assumed to be of finite type over
an algebraically closed field~$\kk$.
For an algebraic variety $X$, we denote by $\D^b(X)$
the bounded derived category of coherent sheaves on $X$,
and by $\D^-(X)$ the unbounded from below derived category
of coherent sheaves.
For $F,G \in \D^-(X)$, we denote by $\RCHom(F,G)$ the local $\RCHom$-complex
and by $F \otimes G$ the derived tensor product.
Similarly, for a map $f:X \to Y$, we denote by
$f_*$ the derived pushforward functor and by $f^*$ the derived pullback.
Finally, $f^!$ stands for the twisted pullback functor.

\subsection{Semiorthogonal decompositions}

If $\CA$ is a full subcategory of $\CT$ then the {\sf right orthogonal}\/ to $\CA$ in $\CT$
(resp.\ the {\sf left orthogonal}\/ to $\CA$ in $\CT$) is the full subcategory
$\CA^\perp$ (resp.\ ${}^\perp\CA$) consisting of all objects $T\in\CT$
such that $\Hom_\CT(A,T) = 0$ (resp.\ $\Hom_\CT(T,A) = 0$) for all $A \in \CA$.

\begin{definition}[\cite{BK,BO1,BO2}]
A sequence $\CA_1,\dots,\CA_n$ of full triangulated subcategories in a triangulated
category $\CT$ is called a {\sf semiorthogonal collection}\/ if $\Hom_{\CT}(\CA_i,\CA_j) = 0$ for $i > j$.
A semiorthogonal collection $\CA_1,\dots,\CA_n$ is a {\sf semiorthogonal decomposition}\/ of $\CT$
if for every object $T \in \CT$ there exists a chain of morphisms
$0 = T_n \to T_{n-1} \to \dots \to T_1 \to T_0 = T$ such that
the cone of the morphism $T_k \to T_{k-1}$ is contained in $\CA_k$
for each $k=1,2,\dots,n$. In other words, if there exists a diagram
\begin{equation}\label{tower}
\vcenter{
\xymatrix@C-7pt{
0 \ar@{=}[r] & T_n \ar[rr]&& T_{n-1} \ar[dl] \ar[rr]&& \quad\dots\quad \ar[rr]&& T_2 \ar[rr]&& T_1 \ar[dl] \ar[rr]&& T_0 \ar[dl] \ar@{=}[r] & T \\
&& A_n \ar@{..>}[ul] &&& \dots &&& A_2 \ar@{..>}[ul]&& A_1 \ar@{..>}[ul]&&
}}
\end{equation}
where all triangles are distinguished (dashed arrows have degree $1$) and $A_k \in \CA_k$.
\end{definition}

Thus, every object $T\in\CT$ admits a decreasing ``filtration''
with factors in $\CA_1$, \dots, $\CA_n$ respectively.
Semiorthogonality implies that this filtration is unique
and functorial. We denote by
$\alpha_k:\CT \to \CA_k$ the functor $T \mapsto A_k$.
We call $\alpha_k$ the \emph{$k$-th projection functor} of the semiorthogonal decomposition.

For any sequence of subcategories $\CA_1,\dots,\CA_n$ in $\CT$ we denote
by $\lan\CA_1,\dots,\CA_n\ran$ the minimal triangulated subcategory of $\CT$
containing $\CA_1$, \dots, $\CA_n$, and by $\lan\CA_1,\dots,\CA_n\ran^\oplus$
the minimal Karoubian (i.e. closed under direct summands) triangulated
subcategory of $\CT$ containing $\CA_1$, \dots, $\CA_n$.

If $\CT = \lan\CA_1,\dots,\CA_n\ran$ is a semiorthogonal decomposition then
$\CA_i = {}^\perp\lan\CA_1,\dots,\CA_{i-1}\ran \cap \lan\CA_{i+1},\dots,\CA_n\ran{}^\perp$.

\begin{definition}[\cite{BK,B}]
A full triangulated subcategory $\CA$ of a triangulated category $\CT$ is called
{\sf right admissible}\/ if for the inclusion functor $i:\CA \to \CT$ there is
a right adjoint $i^!:\CT \to \CA$, and
{\sf left admissible}\/ if there is a left adjoint $i^*:\CT \to \CA$.
Subcategory $\CA$ is called {\sf admissible}\/ if it is both right and left admissible.
\end{definition}

\begin{lemma}[\cite{B}]\label{sod_adm}
If $\CT = \lan\CA,\CB\ran$ is a semiorthogonal decomposition then
$\CA$ is left amissible and $\CB$ is right admissible.
If\/ $\CA_1,\dots,\CA_n$ is a semiorthogonal sequence in $\CT$
such that $\CA_1,\dots,\CA_k$ are left admissible and
$\CA_{k+1},\dots,\CA_n$ are right admissible then
$$
\lan\CA_1,\dots,\CA_k,
{}^\perp\lan\CA_1,\dots,\CA_k\ran \cap \lan\CA_{k+1},\dots,\CA_n\ran{}^\perp,
\CA_{k+1},\dots,\CA_n\ran
$$
is a semiorthogonal decomposition.
\end{lemma}

Let $\CT = \lan \CA,\CB \ran$ be a semiorthogonal decomposition.
Denote by $\alpha:\CA \to \CT$ and $\beta:\CB \to \CT$
the embedding functors. Then for any object $T \in \CT$
we have a distinguished triangle
\begin{equation}\label{sod_tr}
\beta\beta^! T \to T \to \alpha\alpha^*T.
\end{equation}
In particular, $\alpha\alpha^*T$ and $\beta\beta^!T$ are the components
of $T$ with respect to the semiorthogonal decomposition.
So, in this case the projection functors of the semiorthogonal
decompositions are the adjoint functors $\alpha^!$ and $\beta^*$ respectively.

Let $f:X \to S$ be a morphism of algebraic varieties.
Trianguated subcategory $\CA \subset \D^b(X)$ is called \emph{$S$-linear} \cite{K1}
if it is stable with respect to tensoring by pull-backs of perfect complexes on $S$,
i.e.\ if $\CA \otimes f^*(\D^\perf(S)) \subset \CA$.

Let $f:X \to S$ be a projective morphism.
A vector bundle $E$ on $X$ is said to be {\em exceptional over $S$}\/
if it has finite $\Tor$-dimension over $S$ (see~\cite{K1})
and $f_*(E^*\otimes E) \cong \CO_S$.
Note that if on $X$ there exists a bundle of finite $\Tor$-dimension over $S$
then $X$ itself has finite $\Tor$-dimension over $S$.

\begin{lemma}\label{excrel}
If $E$ is a vector bundle on $X$ exceptional over $S$ then the functor
$f_E^*:\D^b(S) \to \D^b(X)$, $G \mapsto f^*G\otimes E$ is fully faithful
and the image of $f_E^*$ is an admissible subcategory in $\D^b(X)$.
\end{lemma}
\begin{proof}
Note that for any $F\in\D^b(X)$, $G\in\D^b(S)$ we have
$$
\begin{array}{l}
\Hom(f^*G\otimes E,F) =
\Hom(f^*G,E^*\otimes F) =
\Hom(G,f_*(E^*\otimes F)),\\[2pt]
\Hom(F,f^*G\otimes E) =
\Hom(F\otimes E^*,f^*G) =
\Hom(F\otimes E^*\otimes f^!\CO_S,f^!G) =
\Hom(f_*(F\otimes E^*\otimes f^!\CO_S),G)
\end{array}
$$
(the second line uses finiteness of $\Tor$-dimension of $f$), so the functor $f_E^*$
has both left and right adjoints. Finally, we have
$$
f_*(E^*\otimes (f^*G\otimes E)) \cong
f_*((E^*\otimes E) \otimes f^*G) \cong
f_*(E^*\otimes E) \otimes G \cong
G
$$
since $E$ is exceptional over $S$. Hence $f_E^*$ is fully faithful.
\end{proof}

\subsection{T-structures and cohomological amplitude of a functor}

Recall that a t-structure on a triangulated category is a pair of full subcategories
$(\CT^{\le 0},\CT^{\ge 0})$ such that
\begin{itemize}
\item $\CT^{\le 0}[1] \subset \CT^{\le 0}$,
$\CT^{\ge 0}[-1] \subset \CT^{\ge 0}$;
\item $\Hom(\CT^{\le 0},\CT^{\ge 1}) = 0$;
\item for any object $T \in \CT$ there exists a distinguished triangle
$$
T^{\le 0} \to T \to T^{\ge 1}
$$
with $T^{\le 0} \in \CT^{\le 0}$ and $T^{\ge 1} \in \CT^{\ge 1}$.
\end{itemize}
As usually, we denote $\CT^{\le k} = \CT^{\le 0}[-k]$, $\CT^{\ge k} = \CT^{\ge 0}[-k]$, $\CT^{[k,l]} = \CT^{\ge k} \cap \CT^{\le l}$.
A t-structure $(\CT^{\le 0},\CT^{\ge 0})$ is \emph{bounded} if
$\CT = \cup_{k,l\in\ZZ} \CT^{[k,l]}$.

Associating to an object $T\in\CT$ the components $T^{\le 0}$ and $T^{\ge 1}$
of $T$ in subcategories $\CT^{\le 0}$ and $\CT^{\ge 1}$ respectively,
one obtaines functors $\tau^{\le 0}:\CT \to \CT^{\le 0}$,
$\tau^{\ge 1}:\CT \to \CT^{\ge 1}$ called the \emph{truncation functors}.
Similarly, we have functors
$\tau^{\le k}(T) = (\tau^{\le 0}(T[k]))[-k]$,
$\tau^{\ge k}(T) = (\tau^{\ge 0}(T[k]))[-k]$,
and the cohomology functors
$\CH^k(T) = \tau^{\ge 0}\tau^{\le 0}(T[k]) \cong \tau^{\le 0}\tau^{\ge 0}(T[k])$.

Let $\CT = \D(\CC)$ be the (bounded or unbounded) derived category of an abelian category $\CC$.
Denote by $\D^{\le 0}(\CC)$ the full subcategory of $\D(\CC)$ formed by all objects
with trivial cohomology in positive degrees, and by $\D^{\ge 0}(\CC)$ the full subcategory
of $\D(\CC)$ formed by all objects with trivial cohomology in negative degrees.
Then $(\D^{\le 0}(\CC),\D^{\ge 0}(\CC))$ is a t-structure, called the \emph{standard t-structure}.
The truncation functors of the standard t-structure are given by the canonical truncations
of complexes, and the cohomology functors by the usual cohomology of complexes.

Let $\Phi:\CS \to \CT$ be a triangulated functor between triangulated categories
equipped with t-structures $(\CS^{\le 0},\CS^{\ge 0})$ and $(\CT^{\le 0},\CT^{\ge 0})$.
The functor $\Phi$ is called \emph{left-exact} if $\Phi(\CS^{\ge 0}) \subset \CT^{\ge 0}$
and \emph{right-exact} if $\Phi(\CS^{\le 0}) \subset \CT^{\le 0}$. More generally,
we will say that $\Phi$ has \emph{finite cohomological amplitude} if there exist
integers $a,b$ such that
$$
\Phi(\CS^{\ge 0}) \subset \CT^{\ge a}
\qquad\text{and}\qquad
\Phi(\CS^{\le 0}) \subset \CT^{\le b}.
$$

\begin{proposition}\label{caf}
Let $X$ be a smooth quasiprojective variety and $\CT$ a triangulated
category with a bounded t-structure $(\CT^{\le 0},\CT^{\ge 0})$.
Then every triangulated functor $\Phi:\D^b(X) \to \CT$
has finite cohomological amplitude.
\end{proposition}
\begin{proof}
Let $X \to \PP^N$ be a locally closed embedding.
Let $\CO_X(k)$ denote the restriction of the line bundle $\CO(k)$ from $\PP^N$ to $X$.
Since the t-structure on $\CT$ is bounded for every $k$
there exist integers $a_k,b_k\in\ZZ$ such that $\Phi(\CO_X(k)) \in \CT^{[a_k,b_k]}$.
Let $a' = \min\{a_0,a_1,\dots,a_N\}$, $b' = \max\{b_0,b_1,\dots,b_N\}$,
so that we have $\Phi(\CO_X(k)) \in \CT^{[a',b']}$ for $0\le k\le N$.
Note that for every $k > N$ there exists a resolution
$$
0 \to \CO^{\oplus c_{k0}} \to \CO(1)^{\oplus c_{k1}} \to \dots \to \CO(N)^{\oplus c_{kN}} \to \CO(k) \to 0
$$
and for every $k < 0$ there exists a resolution
$$
0 \to \CO(k) \to \CO^{\oplus c_{k0}} \to \CO(1)^{\oplus c_{k1}} \to \dots \to \CO(N)^{\oplus c_{kN}} \to 0
$$
on $\PP^N$. Restricting these resolutions to $X$ we deduce that
$\Phi(\CO_X(k)) \in \CT^{[a'-N,b'+N]}$ for all $k\in\ZZ$.
Now let $F$ be any object in $\D^b(X)$.
If $F \in \D^{[p,q]}(X)$ then $F$ admits a locally free left resolution
$$
\{ \dots \to P^{q-2} \to P^{q-1} \to P^q \} \cong F,
$$
where $P^t = \CO_X(k_t)^{c_t}$ sits in degree $t$.
Let $P' = \{P^{p-\dim X} \to P^{p-\dim X + 1} \to \dots \to P^{q-1} \to P^q\}$
be the complex $P^\bullet$ truncated at degree $p - \dim X$.
Then we have a distinguished triangle
$$
P' \to F \to \CH^{p-\dim X}(P')[\dim X - p + 1].
$$
Since $F \in \D^{\ge p}(X)$ and $\CH^{p-\dim X}(P')[\dim X - p + 1] \in \D^{\le p - \dim X - 1}(X)$
and $X$ is smooth we deduce that the map
$F \to \CH^{p-\dim X}(P')[\dim X - p + 1]$ vanishes,
hence $F$ is a direct summand of $P'$.
Finally, we see that $\Phi(P') \in \CT^{[a'-N+p-\dim X,b'+N+q]}$,
hence $\Phi(F) \in \CT^{[a'-N+p-\dim X,b'+N+q]}$ as well.
Therefore $\Phi(\D^{\le 0}(X)) \subset \CT^{\le b'+N}$
and $\Phi(\D^{\ge 0}(X)) \subset \CT^{\ge a'-N-\dim X}$,
so $\Phi$ has finite cohomological amplitude.
\end{proof}

\subsection{Perfect complexes and $\Ext$-amplitude}

Recall that an object $F \in \D(Y)$ in the derived category
of coherent sheaves on an algebraic variety $Y$ is said to be
a \emph{perfect complex} if it is locally (in the Zariski topology)
quasiisomorphic to a bounded complex of locally free sheaves of finite rank.

A similar definition can be given for the category $\D(Y,\CA)$,
the bounded derived category of coherent $\CA$-modules on $Y$,
where $\CA$ is a coherent sheaf of $\CO_Y$-algebras on $Y$.
An object $F \in \D(Y,\CA)$ is said to be a \emph{perfect complex}
if it is locally quasiisomorphic to a bounded complex of locally projective
$\CA$-modules of finite rank.

The category of perfect complexes in $\D(Y,\CA)$ is denoted by $\D^\perf(Y,\CA)$.
It is a triangulated subcategory of $\D(Y,\CA)$. Perfect complexes in $\D^b(Y)$
can be characterized intrinsically by the following property.

\begin{lemma}[\cite{O2}, Pr.~1.11]\label{perfectdby}
An object $F \in \D^b(Y)$ is a perfect complex if and only if for any $G\in\D^b(Y)$
there exists only finite number of $n\in\ZZ$ such that $\Hom(F,G[n])\ne 0$.
\end{lemma}

We will also need a characterization of perfect complexes in the unbounded from below
derived category $\D^-(Y,\CA)$. For this we need the following

\begin{definition}[cf.~\cite{K1}]
An object $F \in \CT$ in a triangulated category $\CT$ with a t-structure $(\CT^{\le 0},\CT^{\ge 0})$
\emph{has finite $\Ext$-amplitude} if there exists an integer $a$ such that
for all $n\in\ZZ$ and all $T \in \CT^{\le n}$ we have
$\RHom(F,T) \in \D^{\le n+a}(\kk)$.
\end{definition}

It is clear that every perfect complex has finite $\Ext$-amplitude
with respect to the standard t-structure.
Now we shall prove the converse.

\begin{lemma}\label{hub_loc}
Let\/ $Y$ be a quasiprojective algebraic variety,
let\/ $\CA$ be a sheaf of coherent $\CO_Y$-algebras, and
let\/ $\CT = \D^-(Y,\CA)$ or $\CT = \D^b(Y,\CA)$ with the standard t-structure.
An object $F \in \CT$ has finite $\Ext$-amplitude if and only if
there exists $a\in\ZZ$ such that
for all $n\in\ZZ$ and all $T \in \CT^{\le n}$ we have
$\RCHom_\CA(F,T) \in \D^{\le n+a}(Y)$.
\end{lemma}
\begin{proof}
Note that
$$
\RHom_\CA(F,T) = \RG(Y,\RCHom_\CA(F,T)).
$$
So, if $\RCHom_\CA(F,T) \in \D^{\le n+a}(Y)$ then
$\RHom(F,T) \in \D^{\le n+a+\dim Y}(\kk)$.

Conversely, assume that $\RHom_\CA(F,T) \in \D^{\le n+a}(\kk)$
for all $n\in\ZZ$ and $T\in\CT^{\le n}$ and let us show that
$\RCHom_\CA(F,T) \in \D^{\le n+a}(Y)$.
Indeed, assume that the complex $C = \RCHom_\CA(F,T)$
has a non-trivial cohomology in degree $m > n+a$.
It is clear that we can choose $k\in\ZZ$ such that
$H^0(Y,\CH^m(C(k))) \ne 0$ and $H^{>0}(Y,\CH^t(C(k))) = 0$
for $t=m-1,m-2,\dots,m-\dim Y + 1$. Then the hypercohomology
spectral sequence of $C(k)$ shows that
$0 \ne H^m(Y,C(k)) = \CH^m(\RHom_\CA(F,T(k)))$, so
$\RHom_\CA(F,T(k)) \not\in \D^{\le n+a}(\kk)$ while $T(k) \in \CT^{\le n}$.
\end{proof}

\begin{proposition}\label{perfect}
Let\/ $Y$ be a quasiprojective algebraic variety,
and let\/ $\CA$ a coherent sheaf of $\CO_Y$-algebras.
Then every bounded object of a finite $\Ext$-amplitude
in $\D^b(Y,\CA)$ or in $\D^-(Y,\CA)$ is a perfect complex.
\end{proposition}
\begin{proof}
Let $F \in \D^-(Y,\CA)$ be a bounded object of finite $\Ext$-amplitude.
Let $a$ be an integer such that
$\RCHom_\CA(F,T) \in \D^{\le n+a}(Y)$ for all $n\in \ZZ$ and all $T \in \D^{\le n}(Y,\CA)$.

Let $P^\bullet = \{\dots \to P^{p-1} \to P^p\}$ be a bounded above complex
of locally free $\CA$-modules of finite rank quasiisomorphic to $F$.
Take $k$ such that $\tau^{\le k}(F) = 0$ and $k \le -a$.
Let $P' = \{P^{k} \to P^{k + 1} \to \dots \to P^{p-1} \to P^p\}$
be the complex $P^\bullet$ truncated at degree $k$.
Let $K = \CH^k(P')$.
Then we have a distinguished triangle
$$
K \to P'[k] \to F[k].
$$
Applying to it the functor $\RCHom_\CA(-,T)$ we obtain a distinguished triangle
$$
\RCHom_\CA(F,T)[-k] \to \RCHom_\CA(P'[k]),T) \to \RCHom_\CA(K,T).
$$
Since $\RCHom_\CA(F,T)[-k] \in \D^{\le n+a+k}(Y) \subset \D^{\le n}(Y)$
and $\RCHom_\CA(P'[k]),T) \le \D^{\le n}(Y)$ we deduce
that $\RCHom_\CA(K,T) \in \D^{\le n}(Y)$ for all $T \in \D^{\le n}(Y,\CA)$.
In particular, $\CExt^{> 0}_\CA(K,T) = 0$ for any sheaf of $\CA$-modules $T$,
hence $K$ is locally projective.
\end{proof}

We also will need the following properties.

\begin{lemma}\label{perforth}
If $X$ is a quasiprojective variety and for $F\in\D(X)$ we have $\Hom(P,F) = 0$
for all $P \in \D^\perf(X)$, then $F = 0$.
\end{lemma}
\begin{proof}
Assume that $\CH^k(F) \ne 0$. Choose $n\in\ZZ$ in such a way
that $H^0(X,\CH^k(F)(n)) \ne 0$ and $H^{>0}(X,\CH^l(F)(n)) = 0$ for $k-\dim X < l < k$.
Then it is easy to see that $H^0(X,\CH^k(F)(n))$ survives in the hypercohomology
spectral sequence, hence $0 \ne H^k(X,F(n)) = \Hom(\CO_X(-n)[-k],F)$.
\end{proof}

\begin{lemma}\label{dminus}
Let $X$ be a quasiprojective variety.
If $G$ is a complex on $X$ unbounded from below, that is
$G \in \D^-(X) \setminus \D^b(X)$, then there exists a line bundle $L$ on $X$
such that $\RHom(L,X)$ is unbounded from below.
\end{lemma}
\begin{proof}
Let $j:X \to X'$ be an open embedding of $X$ into a projective variety $X'$
and let $i:X' \to \PP^N$ be a closed embedding of $X'$ into a projective space.
Then $i_*j_*G$ is an unbounded from below complex of quasicoherent sheaves on $\PP^N$.
Let us check that there exists a line bundle $\CO(t)$ on $\PP^N$ such that
$\RHom(\CO(t),i_*j_*G)$ is unbounded from below. Then by adjunction
$\RHom(j^*i^*\CO(t),G) = \RHom(\CO(t),i_*j_*G)$ is unbounded from below.

For this, consider the Beilinson spectral sequence for $i_*j_*G$
$$
E_1^{p,q} = H^q(\PP^N,i_*j_*G \otimes \CO(-p)) \otimes \Omega^p(p)
\Longrightarrow i_*j_*G.
$$
If for all $0 \le p \le N$ the hypercohomology $H^q(\PP^N,i_*j_*G \otimes \CO(-p))$
vanishes for $q\ll0$ then the spectral sequence implies that $i_*j_*G$ is bounded.
It remains to note that
$H^\bullet(\PP^N,i_*j_*G \otimes \CO(-p)) = \RHom(\CO(p),i_*j_*G)$.
\end{proof}

\subsection{Negative completion}

Let $\CT$ be a triangulated category with a t-structure $(\CT^{\le 0},\CT^{\ge 0})$.
We say that a direct system $T^1 \to T^2 \to T^3 \to \dots $ of objects in $\CT$
\emph{stabilizes in finite degrees}
if for any $k\in\ZZ$ there exists $n_0\in\ZZ$ such that for all $n\ge n_0$
the map $\tau^{\ge k}(T^n) \to \tau^{\ge k}(T^{n+1})$ is an isomorphism.

Let $X$ be an algebraic variety.
Note that any direct system of objects in $\D^b(X)$
which stabilizes in finite degrees with respect to the standard t-structure
has a direct limit in $\D^-(X)$, the unbounded from below derived category.
Moreover, every object in $\D^-(X)$ can be represented as a direct limit
of a stabilizing in finite degrees direct system of perfect complexes
(e.g.\ such system can be obtained by stupid truncations of a locally
free resolution of the object).

\begin{lemma}\label{dminsod}
Assume that $X$ is an algebraic variety and $\D^b(X) = \langle \CA_1,\CA_2,\dots\CA_n\rangle$
is a semiorthogonal decomposition, such that the projection functors $\alpha_k:\D^b(X) \to \CA_k$
have finite cohomological amplitude with respect to the standard t-structure on $\D^b(X)$.
Let $\CA_k^-$ denote minimal triangulated subcategory of $\D^-(X)$ closed under countable
direct sums and containing $\CA_k$.
Then $\D^-(X) = \langle \CA_1^-,\CA_2^-,\dots,\CA_n^-\rangle$
is a semiorthogonal decomposition.
\end{lemma}

Note that the cohomological amplitude of the projection functors $\alpha_k$ is
always finite if $X$ is quasiprojective by proposition~\ref{caf}.

\begin{proof}
First of all, the initial semiorthogonal decomposition induces a semiorthogonal
decomposition of the category of perfect complexes,
$\D^\perf(X) = \langle \CA_1^\perf,\CA_2^\perf,\dots\CA_n^\perf\rangle$
(see~\cite{O2} and \cite{K2}). Now we define
$\CA_k^-$ to be the subcategory of $\D^-(X)$ obtained by iterated addition of cones
to the closure of $\CA_k^\perf$ in $\D^-(X)$ under countable direct sums.
Note that a cone of a direct sum of morphisms is a direct sum of cones of these morphisms,
hence $\CA_k^-$ is a closed under direct sums triangulated dubcategory of $\D^-(X)$.

Now let us check that the categories $\CA_k^-$ form a semiorthogonal decomposition of $\D^-(X)$.
First of all, if $k > l$,
$A_{ki} \in \CA_k^\perf$, $A_{li} \in \CA_l^\perf$ and the direct sums
$\oplus A_{ki}$, $\oplus A_{li}$ exist in $\D^-(X)$ then
$$
\Hom(\oplus_i A_{ki},\oplus_j A_{lj}) =
\Hom(\oplus_i A_{ki},\oplus_j A_{lj}) =
\prod_i \Hom(A_{ki},\oplus_j A_{lj}) =
\prod_i \bigoplus_j \Hom(A_{ki},A_{lj}) = 0
$$
since every perfect complex is a compact object of $\D^-(X)$
(see e.g.~\cite{BV}) and the functor $\Hom(A,-)$ for a compact object $A$
commutes with direct sums.
Addition of cones doesn't spoil semiorthogonality, hence
the collection $\CA_1^-,\CA_2^-,\dots,\CA_n^-$ is semiorthogonal.

It remains to check that the subcategories $\CA_k^-$ generate $\D^-(X)$.
For this we note that $\CA_k^-$ is closed under countable direct limits
(which exist in $\D^-(X)$) since a direct limit can be represented
as a cone of a morphism of direct sums via the telescope construction,
and that every object in $\D^-(X)$ can be represented as a direct limit
of objects in $\D^\perf(X)$. So, take any $T\in\D^-(X)$ and put
$T = \lim\limits_\to T^i$, where $\{T^i\}$ is a stabilizing in finite degrees
direct system in $\D^\perf(X)$.
Let $A_k^i = \alpha_k(T^i)$. It is a direct system in $\CA_k^\perf$.
Since the projection functors $\alpha_k$ have finite cohomological amplitude
it follows that each of these direct systems stabilizes in finite degrees.
Put $A_k = \lim\limits_\to A_k^i$. Then it is easy to see that $A_k \in \CA_k^-$
and that $A_k$ are the projections of $T$ to $\CA_k^-$.
\end{proof}

The categories $\CA_k^-$ constructed in the above lemma are called
\emph{negative completions} of $\CA_k$. It is clear by definition
that $\CA_k^- \supset \CA_k^\perf$. In particular, if $X$ is smooth
then $\CA_k^- \supset \CA_k$. However, for singular $X$ this inclusion
is not clear.

We will need below the following properties of negative completions.

\begin{lemma}\label{dbdm}
Let $\D^-(X) = \langle \CA_1^-,\CA_2^-,\dots,\CA_n^-\rangle$ be the semiorthogonal decomposition
formed by the negative completion of the components of a semiorthogonal decomposition
$\D^b(X) = \langle \CA_1,\CA_2,\dots\CA_n\rangle$. If $X$ is smooth then for any bounded object
$F \in \D^b(X) \subset \D^-(X)$ all the components of $F$ with respect to the decomposition
$\D^-(X) = \langle \CA_1^-,\CA_2^-,\dots,\CA_n^-\rangle$ are bounded.
\end{lemma}
\begin{proof}
We just decompose $F$ with respect to the decomposition $\D^b(X) = \langle \CA_1,\CA_2,\dots\CA_n\rangle$
and note that this gives us a decomposition with respect to $\D^-(X) = \langle \CA_1^-,\CA_2^-,\dots,\CA_n^-\rangle$
since $\CA_k = \CA_k^\perf \subset \CA_k^-$.
\end{proof}


\subsection{Lefschetz decompositions}

Let $X$ be an algebraic variety and $\CO_X(1)$ a line bundle on $X$.

\begin{definition}[\cite{K2}]
A \emph{Lefschetz decomposition} of $\D^b(X)$ is a semiorthogonal decomposition
of the form
$$
\D^b(X) = \lan \CB_0,\CB_1(1),\dots,\CB_{m-1}(m-1)\ran,
\qquad\text{where $0 \subset \CB_{m-1} \subset \dots \subset \CB_1 \subset \CB_0 \subset \D^b(X)$.}
$$
A Lefschet decomposition is called \emph{rectangular}
if all its components coincide
$$
\CB_0 = \CB_1 = \dots = \CB_{m-1}.
$$
Similarly,
a \emph{dual Lefschetz decomposition} of $\D^b(X)$ is a semiorthogonal decomposition
of the form
$$
\D^b(X) = \lan \CB_{m-1}(1-m),\dots,\CB_1(-1),\CB_0\ran,
\qquad\text{where $0 \subset \CB_{m-1} \subset \dots \subset \CB_1 \subset \CB_0 \subset \D^b(X)$.}
$$
\end{definition}

Actually, these notions are equivalent.
Given a Lefschetz decomposition
one can canonically construct a dual Lefschetz decomposition with
the same category $\CB_0$ (as it is shown in the following lemma)
and vice versa.

\begin{lemma}\label{dual_ld}
Assume that
$$
\D^b(X) = \langle \CB_0',\CB_1'(1),\dots,\CB_{m-1}'(m-1) \rangle,
$$
is a Lefschetz decomposition and that all subcategories
$\CB'_0,\CB'_1,\dots,\CB'_{m-1}$ are admissible in $\D^b(X)$. Put
$$
\CB_0 = \CB'_0,
\qquad\text{and}\qquad
\CB_{k} = \CB'_0(k)^\perp \cap \CB_{k-1}
\qquad\text{for $k=1,2,\dots,m-1$}.
$$
Then the chain of subcategories $0 \subset \CB_{m-1} \subset \dots \subset \CB_1 \subset \CB_0$
gives a dual Lefschetz decomposition
$$
\D^b(X) = \langle \CB_{m-1}(1-m),\dots,\CB_1(-1),\CB_0 \rangle.
$$
\end{lemma}
\begin{proof}
For a start let us note that
\begin{equation}\label{bprime}
\langle\CB'_0,\CB'_1(1),\dots,\CB'_k(k)\rangle =
\langle\CB'_0,\dots,\CB'_0(k-1),\CB'_0(k)\rangle.
\end{equation}
Indeed, ``$\subset$'' is evident and ``$\supset$'' follows from
$\Hom(\CB'_p(p),\CB'_0(q)) = \Hom(\CB'_p(p-q),\CB'_0) = 0$
for all $q \le k < p$ since $\CB'_p \subset \CB'_{p-q}$
because $\D^b(X) = \langle \CB_0',\CB_1'(1),\dots,\CB_{m-1}'(m-1) \rangle$
is a semiorthogonal decomposition.

We are going to check by induction that
$$
\langle \CB'_0,\CB'_1(1),\dots,\CB'_k(k) \rangle = \langle \CB_k,\dots,\CB_1(k-1),\CB_0(k)\rangle
$$
(both sides are semiorthogonal decompositions).
The base of the induction, $k=0$, is evident by definition of $\CB_0$.
Assume that the above decomposition is true for $k-1$.
Note that
the collection $\langle\CB_{k-1}(1),\dots,\CB_1(k-1),\CB_0(k)\rangle$
is semiorthogonal by the induction hypothesis and being equivalent to
$\langle\CB'_0(1),\CB'_1(2),\dots,\CB'_{k-1}(k)\rangle$ it is admissible,
so it remains to note that
\begin{multline*}
\langle\CB_{k-1}(1),\dots,\CB_1(k-1),\CB_0(k)\rangle^\perp_{\langle\CB'_0,\CB'_1(1),\dots,\CB'_k(k)\rangle} =
\langle\CB'_0(1),\CB'_1(2),\dots,\CB'_{k-1}(k)\rangle^\perp_{\langle\CB'_0,\CB'_1(1),\dots,\CB'_k(k)\rangle} =
\\ =
\langle\CB'_0(1),\CB'_0(2),\dots,\CB'_0(k)\rangle^\perp_{\langle\CB'_0,\CB'_0(1),\dots,\CB'_0(k)\rangle} =
\langle\CB'_0(1),\CB'_0(2),\dots,\CB'_0(k-1)\rangle^\perp_{\langle\CB'_0,\CB'_0(1),\dots,\CB'_0(k-1)\rangle} \cap \CB'_0(k)^\perp
\end{multline*}
(the first equality is by the induction hypothesis, the second equality is by~\eqref{bprime}),
and that
\begin{multline*}
\langle\CB'_0(1),\CB'_0(2),\dots,\CB'_0(k-1)\rangle^\perp_{\langle\CB'_0,\CB'_0(1),\dots,\CB'_0(k-1)\rangle} =
\langle\CB'_0(1),\CB'_1(2),\dots,\CB'_{k-2}(k-1)\rangle^\perp_{\langle\CB'_0,\CB'_1(1),\dots,\CB'_{k-1}(k-1)\rangle} =
\\ =
\langle\CB_{k-2}(1),\dots,\CB_1(k-2),\CB_0(k-1)\rangle^\perp_{\langle\CB'_0,\CB'_1(1),\dots,\CB'_{k-1}(k-1)\rangle} =
\CB_{k-1}
\end{multline*}
(the first equality is by~\eqref{bprime}, the second and the third equalities are by the induction hypothesis).
\end{proof}

\section{Categorical resolutions of singularities}\label{defres}

To define a resolution of singularities on a categorical level one should
first define a class of triangulated categories to be considered
as analogs of derived categories of smooth varieties.
Presently, there are several approaches see \cite{BO2,TV,KS}.
Unfortunately, \cite{BO2} seams to be too restrictive,
while \cite{TV} and \cite{KS} are not sufficiently worked out.
So, for the purpose of this paper we use the following definition.

\begin{definition}
A triangulated category $\D$ is {\em regular} if it is equivalent
to an admissible subcategory of the bounded derived category of a
smooth algebraic variety.
\end{definition}

Certainly, this definition should be considered as provisional.

One more ingredient we need is a notion of a perfect object
in the resolved triangulated category. Here we use the property
of perfect complexes stated in lemma~\ref{perfectdby}.
We say that an object $F\in\D$ is a perfect complex if for any $G\in\D$
there exists only finite number of $n\in\ZZ$ such that $\Hom(F,G[n])\ne 0$.

\begin{definition}
A \emph{categorical resolution} of a triangulated category $\D$
is a regular triangulated category $\TD$ and a pair of functors
$$
\pi_*:\TD \to \D,
\qquad
\pi^*:\D^\perf \to \TD,
$$
such that $\pi^*$ is left adjoint to $\pi_*$ on $\D^\perf$, that is
$$
\Hom_\TD(\pi^*F,G) \cong \Hom_\D(F,\pi_*G)
\qquad\text{for any $F\in\D^\perf$, $G\in\TD$},
$$
and the natural morphism of functors $\id_{\D^\perf} \to \pi_*\pi^*$
is an isomorphism.
\end{definition}

\begin{remark}
The last assumption restricts the class of singularities in consideration
to rational singularities. However, in general one needs a kind of
surjectivity assumption for the functor $\pi_*$ (or injectivity for $\pi^*$)
to be made (otherwise, e.g. $\TD = 0$ will satisfy the definition).
And even in the case of usual resolution of singularities, it is not clear
what kind of surjectivity is satisfied.
\end{remark}

Let $\D = \D^b(Y)$ be the bounded derived category of an algebraic variety $Y$
with rational singularities.
If $\pi:\TY \to Y$ is a resolution of singularities of $Y$ then
the category $\TD = \D^b(\TY)$ with the pushforward $\pi_*:\D^b(\TY) \to \D^b(Y)$
and the pullback $\pi^*:\D^\perf(Y) \to \D^b(\TY)$ functors form a categorical
resolution of $\D^b(Y)$. Note that the pullback on the whole $\D^b(Y)$
takes values in $\D^-(\TY)$. This is why we ask  for the functor $\pi^*$
to be defined only on the subcategory $\D^\perf$.

Another example of a categorical resolution is given by
noncommutative resolutions of singularities introduced by Van den Bergh~\cite{V1,V2}.
If $Y$ is a singular projective algebraic variety and $\CA$ is a sheaf of noncommutative
algebras on $Y$ giving a noncommutative resolution of singularities of $Y$ then
the category $\TD = \D^b(\CA\!-\!\mod)$ is a categorical resolution of $\D^b(Y)$.

\begin{definition}
A categorical resolution $(\TD,\pi_*,\pi^*)$ of $\D$
is \emph{crepant} if the functor $\pi^*$ is right adjoint to $\pi_*$ on $\D^\perf$
$$
\Hom_\TD(G,\pi^*F) \cong \Hom_\D(\pi_*G,F)
\qquad\text{for any $F\in\D^\perf$, $G\in\TD$}.
$$
\end{definition}

Certainly, given a geometrical resolution of singularities $\pi:\TY \to Y$,
the corresponding categorical resolution $(\D^b(\TY),\pi_*,\pi^*)$
is crepant if and only if $\pi$ is crepant.
Similarly, the categorical resolution corresponding to
a noncommutative crepant resolution in the sense of Van den Bergh
is also crepant.


\section{Lefschetz decompositions and categorical resolutions}\label{ld_cr}

The goal of this section is to give an example of a categorical resolution
of the derived category of coherent sheaves on a singular algebraic variety
satisfying some special conditions.

We start with the following observation.

\begin{proposition}\label{isld}
Let $D \subset X$ be a Cartier divisor, let $i:D \to X$ be the embedding, let $L = \CN^*_{D/X}$ be the conormal bundle,
and let $\D^b(D) = \langle \CB_{m-1}\otimes L^{1-m},\CB_{m-2}\otimes L^{2-m},\dots,\CB_1\otimes L^{-1},\CB_0 \rangle$
be a dual Lefschetz decomposition with respect to $L$. Then the pushforward functor $i_*$ is fully faithful on
the subcategories $\CB_k\otimes L^{-k}$ for $1\le k\le m-1$ and we have a semiorthogonal decomposition
$$
D^b(X) = \langle i_*(\CB_{m-1}\otimes L^{1-m}),i_*(\CB_{m-2}\otimes L^{2-m}),\dots,i_*(\CB_1\otimes L^{-1}),\TD \rangle,
$$
where $\TD = \{ F \in \D^b(X)\ |\ i^*F \in \CB_0 \}$.

Similarly, if $\D^b(D) = \langle \CB'_0,\CB'_1\otimes L,\dots,\CB'_{m-2}\otimes L^{m-2},\CB'_{m-1}\otimes L^{m-1} \rangle$
is a Lefschetz decomposition with respect to $L$ then $i_*$ is fully faithful on the subcategories $\CB'_k\otimes L^k$
for $1\le k\le m-1$ and we have a semiorthogonal decomposition
$$
D^b(X) = \langle \TD',i_*(\CB'_1\otimes L),\dots,i_*(\CB'_{m-2}\otimes L^{m-2}),i_*(\CB'_{m-1}\otimes L^{m-1}) \rangle,
$$
where $\TD' = \{ F \in \D^b(X)\ |\ i^!F \in \CB'_0 \}$.

In particular, if all subcategories $\CB_k \subset \D^b(D)$ are admissible then $\TD$ is admissible.
\end{proposition}
\begin{proof}
Consider any objects $F \in \CB_k\otimes L^{-k}$, $G \in \CB_l\otimes L^{-l}$, $1\le k\le l\le m-1$.
Note that $\Hom(i_*F,i_*G) = \Hom(i^*i_*F,G)$ by adjunction.
On the other hand, $i$ is a divisorial embedding and $\CN_{\TZ/\TY} \cong L^{-1}$.
Therefore for any object $F \in \D^b(\TZ)$ we have a distinguished triangle
$$
i^*i_*F \to F \to F\otimes L[2].
$$
Applying $\RHom(-,G)$ functor we obtain a distinguished triangle
$$
\RHom(F\otimes L[2],G) \to \RHom(F,G) \to \RHom(i_*F,i_*G).
$$
But for $F \in \CB_k\otimes L^{-k}$ we have
$$
F\otimes L[2] \in \CB_k\otimes L^{1-k} \subset \CB_{k-1}\otimes L^{1-k}.
$$
Therefore $F\otimes L[2]$ is left orthogonal to $G$ by~(\ref{ldtz}).
We deduce that $\Hom(i_*F,i_*G) = \Hom(F,G)$.
In particular, the functor $i_*$ is fully faithful on
subcategories $\CB_k\otimes L^{-k}$ for $1\le k\le m-1$,
and the collection
\begin{equation}\label{iscb}
i_*(\CB_{m-1}\otimes L^{1-m}),i_*(\CB_{m-2}\otimes L^{2-m}),\dots,i_*(\CB_1\otimes L^{-1})
\end{equation}
is semiorthogonal. Furthermore, if an object $F \in \D^b(X)$ is left orthogonal to all these subcategories
then $0 = \Hom(F,i_*(\CB_k\otimes L^{-k})) = \Hom(i^*F,\CB_k\otimes L^{-k})$ for $1\le k\le m-1$,
therefore $i^*F \in \CB_0$. This shows that the left orthogonal to the collection~\eqref{iscb}
coincides with $\TD$. It remains to note that the functor $i_*$ has a right adjoint $i^!:\D^b(X) \to \D^b(D)$,
hence every subcategory $i_*(\CB_k\otimes L^{-k})$ is left admissible in
${}^\perp\langle i_*(\CB_{m-1}\otimes L^{1-m}),i_*(\CB_{m-2}\otimes L^{2-m}),\dots,i_*(\CB_{k-1}\otimes L^{1-k})\rangle$,
hence the first claim.

The second claim can be proved by the same arguments. For the last claim note that
if all categories $\CB'_k$ are admissible then we can use the construction of lemma~\ref{dual_ld}
to obtain a Lefschetz decomposition
$\D^b(D) = \langle \CB'_0,\CB'_1\otimes L,\dots,\CB'_{m-2}\otimes L^{m-2},\CB'_{m-1}\otimes L^{m-1} \rangle$
with $\CB'_0 = \CB$. Then for the left admissible subcategory $\TD' \subset \D^b(X)$
corresponding to this Lefschetz decomposition we have $\TD' = \TD\otimes\CO_X(-D)$.
Indeed, $i^!F \cong i^*F\otimes\CO_D(D)[-1] \cong i^*(F\otimes\CO_X(D))[-1]$
hence $i^!F \in \CB'_0$ is equivalent to $i^*(F\otimes\CO_X(D)) \in \CB_0$,
so $F \in \TD'$ is equivalent to $F\otimes\CO_X(D) \in \TD$.
\end{proof}


We are going to apply proposition~\ref{isld} for the construction of a categorical
resolution of singularities of $\D^b(Y)$, where $Y$ is an algebraic variety
with rational singularities. We consider a usual resolution $\pi:\TY \to Y$
%
such that the exceptional locus of $\pi$ is an irreducible divisor $\TZ \subset \TY$,
take a dual Lefschetz decomposition
\begin{equation}\label{ldtz}
\D^b(\TZ) = \langle \CB_{m-1}(1-m),\CB_{m-2}(2-m),\dots,\CB_1(-1),\CB_0 \rangle,
\end{equation}
with respect to the line bundle
$$
\CO_\TZ(1) := \CN^*_{\TZ/\TY}
$$
and take the corresponding subcategory
\begin{equation}\label{deftd}
\TD := \{ F\in\D^b(\TY)\ |\ i^*F \in \CB_0\} \subset \D^b(\TY)
\end{equation}
to be the resolving category. Note that by proposition~\ref{isld} we
have a semiorthogonal decomposition
\begin{equation}\label{sodty}
D^b(\TY) = \langle i_*(\CB_{m-1}(1-m)),i_*(\CB_{m-2}(2-m)),\dots,i_*(\CB_1(-1)),\TD \rangle,
\end{equation}
hence the subcategory $\TD \subset \D^b(\TY)$ is right admissible, hence $\TD$
is regular by definition. We denote the restriction of the
pushforward functor $\pi_*:\D^b(\TY) \to \D^b(Y)$ to $\TD$ by the
same letter $\pi_*:\TD \to \D^b(Y)$. Similarly, we have the pullback
$\pi^*:\D^\perf(Y) \to \D^b(\TY)$. Now we need a condition for
$\pi^*(\D^\perf(Y))$ to be contained in $\TD$.

\begin{lemma}\label{pisintscb0}
Let $Z = \pi(\TZ)$ and $p:\TZ \to Z$ be the restriction of $\pi$ to $\TZ$.
If $p^*(\D^\perf(Z)) \subset \CB_0$ then $\pi^*(\D^\perf(Y)) \subset \TD$.
\end{lemma}
\begin{proof}
Let $j:Z \to Y$ be the embedding, so that $j\circ p = \pi\circ i$.
Take any $F \in \D^\perf(Y)$ and note that $i^*\pi^*(F) = p^*j^*(F) \in p^*(\D^\perf(Z)) \subset \CB_0$,
hence $\pi^*(\D^\perf(Y)) \in \TD$.
\end{proof}

Now we have a pair of adjoint functors $\pi^*:\D^\perf(Y) \to \TD$ and $\pi_*:\TD \to \D^b(Y)$.
Moreover, $\pi_*\pi^* \cong \id$ by the projection formula since $Y$ has rational singularities.
Thus we have proved the following

\begin{theorem}\label{th1}
Let $\pi:\TY \to Y$ be a resolution of rational singularities.
Assume that $\TZ \subset \TY$ is the exceptional locus of $\pi$,
$Z = \pi(\TZ)$, and $p = \pi_{|\TZ}:\TZ \to Z$.
If~\eqref{ldtz} is a dual Lefschetz decomposition of $\D^b(\TZ)$
with respect to the line bundle $\CO_{\TZ}(1) = \CN^*_{\TZ/\TY}$
and the category $\CB_0 \subset \D^b(\TZ)$ contains $p^*\D^\perf(Z)$ then
the category $\TD \subset \D^b(\TY)$ is a categorical resolution
of $\D^b(Y)$.
\end{theorem}

Though the tensor structure of the category $\D^\perf(Y)$ doesn't appear in the definition
of a categorical resolution of singularities, it is an interesting question when the resolution
$\TD$ of $\D^b(Y)$ admits a module structure over $\D^\perf(Y)$, compatible with the module structure
of $\D^b(Y)$ via the functor $\pi_*$. A natural way to obtain such a structure is by restricting
the natural module structure of $\D^b(\TY)$, and it is clear that this structure restricts
if and only if $\TD$ is a $Y$-linear subcategory in $\D^b(\TY)$.

\begin{lemma}\label{tdlin}
If
the subcategory $\CB_0 \subset \D^b(\TZ)$ is $Z$-linear then
the subcategory $\TD \subset \D^b(\TY)$ is $Y$-linear.
\end{lemma}
\begin{proof}
If $G \in \D^\perf(Y)$ and $F \in \D^b(\TY)$ then
$i^*(F\otimes \pi^*G) \cong
i^*F \otimes i^*\pi^*G \cong
i^*F \otimes p^*j^*G$.
Now, if $\CB_0$ is $Z$-linear then
$i^*F \in \CB_0$ implies $i^*F\otimes p^*j^*G \in \CB_0$, hence
$F \in \TD$ implies $F\otimes \pi^*G \in \TD$.
\end{proof}

\begin{remark}
It is also easy to prove the inverse statement. For this it suffices to note that for any object
$F \in \CB_0 \subset \D^b(\TZ)$ there exists $F' \in \D^b(\TY)$ such that $i^*F' \cong F$.
\end{remark}

It is interesting that crepancy of the categorical resolution $(\TD,\pi_*,\pi^*)$
of $\D^b(Y)$ can be formulated in terms of the Lefschetz decomposition~\eqref{ldtz}.

\begin{proposition}\label{pi_crep}
If the Lefschetz decomposition~\eqref{ldtz} of $\D^b(\TZ)$ is rectangular
{\rm(}i.e.\ $\CB_0 = \CB_1 = \dots = \CB_{m-1}${\rm)}, $Y$ is Gorenstein
and $K_\TY = \pi^*K_Y + (m-1)\TZ$ then $\TD$ is a crepant categorical resolution of $\D^b(Y)$.
\end{proposition}
\begin{proof}
Note that the projection functor $\TD \to \D^b(Y)$ takes form $G \mapsto \pi_*\delta(G)$,
where $\delta:\TD \to \D^b(\TY)$ is the embedding.
Therefore, the adjoint functors are given by
$F \mapsto \delta^*\pi^*F$ and $F \mapsto \delta^!\pi^!F$ respectively
and we have to construct an isomorphism of functors $\delta^!\pi^! \cong \delta^*\pi^*$.

Take any $F \in \D^\perf(Y)$.
First of all, note that $\pi^*F \in \TD$ by lemma~\ref{pisintscb0},
hence $\delta\delta^*\pi^*F = \pi^*F$.
On the other hand consider $\delta\delta^!\pi^!F$.
Since $\pi^!(F) \cong \pi^*(F) \otimes \omega_{\TY/Y} \cong \pi^*(F)\otimes\CO_\TY((m-1)\TZ)$
we deduce from~\eqref{sod_tr}
that $\delta\delta^!\pi^!F$ is just the $\TD$-component of $\pi^*F\otimes\CO_\TY((m-1)\TZ)$
in the semiorthogonal decomposition~\eqref{sodty}. Thus we should check that it coincides
with $\delta\delta^*\pi^*F \cong \pi^*F$.

For this we consider the short exact sequence
$0 \to \CO_\TY(-\TZ) \to \CO_\TY \to i_*\CO_\TZ \to 0$.
Tensoring it by
$\pi^*F\otimes\CO_Y((m-1)\TZ)$,
$\pi^*F\otimes\CO_Y((m-2)\TZ)$, \dots,
$\pi^*F\otimes\CO_Y(\TZ)$,
we obtain distinguished triangles
$$
\begin{array}{ccccl}
\pi^*F\otimes\CO_\TY((m-2)\TZ) &\to&
\pi^*F\otimes\CO_\TY((m-1)\TZ) &\to&
i_*i^*(\pi^*F\otimes\CO_\TY((m-1)\TZ)),\\
\pi^*F\otimes\CO_\TY((m-3)\TZ) &\to&
\pi^*F\otimes\CO_\TY((m-2)\TZ) &\to&
i_*i^*(\pi^*F\otimes\CO_\TY((m-2)\TZ)),\\
&&\vdots\\
\pi^*F &\to&
\pi^*F\otimes\CO_\TY(\TZ) &\to&
i_*i^*(\pi^*F\otimes\CO_\TY(\TZ))
\end{array}
$$
We can rewrite these triangles in the following diagram
$$
\vcenter{
\xymatrix@C-19.6pt{
\pi^*F \ar[rr] && \pi^*F\otimes\CO_\TY(\TZ) \ar[dl] \ar[rr] && \quad\dots\quad \ar[rr]&& \pi^*F\otimes\CO_\TY((m-2)\TZ) \ar[rr]&& \pi^*F\otimes\CO_\TY((m-1)\TZ) \ar[dl] \\
& \qquad\makebox[0pt]{\hss $i_*i^*(\pi^*F\otimes\CO_\TY(\TZ))$\hss}{ \mathstrut\quad} \ar@{..>}[ul] &&& \dots &&&
\quad\makebox[0pt]{\hss $i_*i^*(\pi^*F\otimes\CO_\TY((m-1)\TZ))$\hss}{\mathstrut\quad} \ar@{..>}[ul]&&
}}
$$
Note that
$i^*(\pi^*F\otimes\CO_\TY(k\TZ)) \cong i^*\pi^*F \otimes \CO_{\TZ/Z}(-k)$
and $i^*\pi^*F = p^*j^*F \in p^*(\D^\perf(Z)) \subset \CB_0$ for any $F\in\D^\perf(Y)$.
Therefore
$i^*(\pi^*F\otimes\CO_\TY(k\TZ)) \in \CB_0(-k) = \CB_k(-k)$
since the collection~\eqref{ldtz} is $Z$-linear and rectangular.
Therefore, the above diagram gives a semiorthognal decomposition of $\pi^*F\otimes\CO_\TY((m-1)\TZ)$
with respect to~\eqref{sodty} and $\pi^*F$ is the component of $\pi^*F\otimes\CO_\TY((m-1)\TZ)$
in $\TD$. So, $\delta\delta^!\pi^!F = \pi^*F$
and we are done.
\end{proof}

The same argument gives an information about the Serre functor of $\TD$ if $Y$ is projective.

\begin{proposition}\label{serre}
Assume that $Y$ is Gorenstein and projective, and $K_\TY = \pi^*K_Y + (m-1)\TZ$.
Then the category $\TD$ has a Serre functor $\FS_\TD$, and for any $F \in \TD$
such that $i^*F \in \CB_{m-1}$ we have
$$
\FS_\TD(F) \cong F\otimes \pi^*\CO_Y(K_Y)[\dim Y].
$$
In particular, if the Lefschetz decomposition~\eqref{ldtz} is rectangular,
then $\FS_\TD(F) \cong F\otimes \pi^*\CO_Y(K_Y)[\dim Y]$ for any $F\in\TD$.
\end{proposition}
\begin{proof}
The map $\pi:\TY \to Y$ is projective by definition, hence $\TY$ is projective.
Therefore, $\D^b(\TY)$ has a Serre functor $\FS_\TY$ which is given by
$\FS_\TY(G) = G \otimes \CO_\TY(K_\TY)[\dim\TY]$.
Since $\TD$ is an admissible subcategory in $\D^b(\TY)$
it follows that $\TD$ also has a Serre functor $\FS_\TD$ and that
$\FS_\TD \cong \beta^!\FS_\TY\beta$.

Now take any $F \in \TD$ such that $i^* F \in \CB_{m-1}$. Then
\begin{multline*}
\FS_\TY\beta(F) \cong
F\otimes\CO_\TY(K_\TY)[\dim\TY] \cong
F\otimes\CO_\TY(\pi^*K_Y + (m-1)\TZ)[\dim Y] \cong
\\ \cong
(F\otimes\pi^*\CO_Y(K_Y)[\dim Y]) \otimes \CO_\TY((m-1)\TZ)[\dim Y],
\end{multline*}
and the same arguments as in proposition~\ref{pi_crep}
show that $\beta^!\FS_\TY\beta(F) \cong F\otimes\pi^*\CO_Y(K_Y)[\dim Y]$.
\end{proof}

Recall \cite{BO2} that a categorical resolution of singularities is called
\emph{minimal} if it can be embedded as an admissible subcategory into any
other categorical resolution of the same singularity. It is conjectured
by Bondal and Orlov in \emph{loc.~cit.}\/ that every singularity
admits a minimal categorical resolution.
We expect that the minimality of the categorical resolution of $\D^b(Y)$
corresponding to a given Lefschetz decomposition of $\D^b(X)$
can be expressed in terms of the Lefschetz decomposition as follows.

Let us say that a dual Lefschetz decomposition
$\D^b(X) = \lan \CB_{m-1}\otimes L^{1-m},\dots,\CB_1\otimes L^{-1},\CB_0 \ran$
is \emph{strict}, if
the whole decomposition can be reconstructed from the subcategory $\CB_0$
as in lemma~\ref{dual_ld}, that is if
\begin{equation}\label{ldstr}
\CB_{k} = (\CB_0\otimes L^k)^\perp \cap \CB_{k-1},
\qquad\text{for $k=1,\dots,m-1$.}
\end{equation}
Any strict Lefschetz decomposition is determined uniquely
by the subcategory $\CB_0$. The inclusion relation on the set
of all subcategories of $\D^b(X)$ gives a partial ordering
on the set of all strict Lefschetz decompositions.
A Lefschetz decomposition is called \emph{minimal},
if it is strict and minimal with respect to this partial ordering.

\begin{conjecture}\label{c1}
A categorical resolution of singularities of\/ $Y$ corresponding
to a Lefschet decomposition of\/ $\D^b(\TZ)$ is minimal if and only if
the Lefschetz decomposition is minimal with respect
to the above partial ordering.
\end{conjecture}

\begin{conjecture}\label{c2}
Any indecomposable crepant categorical resolution of\/ $Y$ is minmal.
In particular, all crepant categorical resolutions of\/ $Y$ are equivalent.
\end{conjecture}

Note that a rectangular Lefschetz decomposition is always minimal.

\begin{conjecture}\label{c3}
For every smooth projective variety $X$ and any ample line bundle $L$ on $X$,
there exists a minimal Lefschetz decomposition of $\D^b(X)$.
\end{conjecture}

The latter conjecture is very likely to be true. Indeed,
it suffices to prove that every decreasing chain of admissible
subcategories $\CB_0 \supset \CB'_0 \supset \dots \supset \CB^{(n)}_0 \supset \dots$
stabilizes.
Conjectures~\ref{c1} and~\ref{c3} together imply the conjecture
of Bondal and Orlov for cones over smooth varieties.
Conjecture~\ref{c3} also is very important from the point of view of theory
of Homological Projective Duality~\cite{K2}.

\section{Noncommutative resolutions}\label{ncr}

In this section we investigate when the categorical resolution $\TD$ of $\D^b(Y)$
constructed in the previous section can be realized as the derived category
of sheaves of modules over a sheaf of $\CO_Y$-algebras on $Y$. In other words,
when $\TD$ is a noncommutative resolution of singularities? Let us formulate
the main result.

As in the previous section, assume that $\pi:\TY \to Y$ is a resolution of rational singularity,
the exceptional locus $\TZ \subset \TY$ of which is an irreducible divisor.
Let $Z = \pi(\TZ) \subset Y$ and denote by $p:\TZ \to Z$ the projection,
and by $i:\TZ \to \TY$, $j:Z \to Y$ the embeddings, so that we have the following
commutative diagram
\begin{equation}\label{yz}
\vcenter{\xymatrix{
\TZ \ar[r]^i \ar[d]_{p} & \TY \ar[d]^\pi \\
Z \ar[r]^j & Y
}}
\end{equation}
Assume that~\eqref{ldtz} is a Lefschetz decomposition of $\D^b(\TZ)$ with respect
to the line bundle $\CO_\TZ(1):=\CN^*_{\TZ/\TY}$, such that $p^*(\D^\perf(Z)) \subset \CB_0$
and all subcategories $\CB_0,\CB_1,\dots,\CB_{m-1}$ are both left and right admissible.
Note that this is always true if each $\CB_k$ is generated by an exceptional over $Z$
collection of vector bundles as in lemma~\ref{excrel}.

\begin{remark}\label{tdadm}
If all subcategories $\CB_0,\CB_1,\dots,\CB_{m-1}$ are both left and right admissible then
by proposition~\ref{isld} the categorical resolution $\TD$ of $\D^b(Y)$ defined by~\eqref{deftd}
is also both left and right admissible subcategory of $\D^b(\TY)$,
and we have a semiorthogonal decomposition~\eqref{sodty}.
\end{remark}

Assume {\em additionally}\/ that the category $\CB_0$ of the Lefschetz decomposition~\eqref{ldtz}
is generated over $\D^\perf(Z)$ by a vector bundle on $\TY$. In other words, assume that there
is a vector bundle $E$ on $\TY$ such that
\begin{equation}\label{ldtz1}
\CB_0 = \langle i^*E \otimes p^*(\D^b(Z)) \rangle^\oplus,
\end{equation}
that is $i^*E$ has finite $\Tor$-dimension over $Z$ and $\CB_0$ is the minimal Karoubian
(i.e.\ closed under direct summands) triangulated subcategory in $\D^b(\TZ)$
containing $i^*E\otimes p^*(\D^b(Z))$.
Note that~\eqref{ldtz1} implies that $\CB_0$ is $Z$-linear, hence $\TD$ is $Y$-linear by lemma~\ref{tdlin}.

Assume also, that $E$ is \emph{tilting over $Y$}, that is
the pushforward $\pi_*\CEnd E$ is a pure sheaf on $Y$.
Then
\begin{equation}\label{defa}
\CA = \pi_*\CEnd E \in \Coh(Y),
\end{equation}
is a sheaf of $\CO_Y$-algebras on $Y$.

Finally, assume that $\TZ$ is a scheme-theoretic preimage of $Z$, that is
\begin{equation}\label{pisjz}
\pi^{-1}\CJ_Z\cdot\CO_\TY = \CJ_\TZ,
\end{equation}
where $\CJ_Z$ and $\CJ_\TZ$ are the ideals of $Z$ in $Y$ and $\TZ$ in $\TY$ respectively.
Then we have the following

\begin{theorem}\label{th2}
Assume that the above conditions are satisfied.
Then the sheaf of algebras $\CA = \pi_*\CEnd E$ has finite homological dimension
and the category $\D^b(Y,\CA)$ is a categorical resolution of\/ $\D^b(Y)$.
\end{theorem}

The proof takes the rest of the section.

Consider the functor
$$
\Phi:\D^-(\TY) \to \D^-(Y,\CA),\qquad
F \mapsto \pi_*(F \otimes E^*).
$$
Below we will show that $\Phi$ induces an equivalence $\TD \to \D^b(Y,\CA)$.
For a start we need a description of an adjoint functor to $\Phi$.
A naive approach would be to consider $\Psi(G) = \pi^*G \otimes_{\pi^*\CA} E$.
However, even to define the functor $\Psi$ in this way is difficult.
The main problem is that $\pi^*\CA$ is a DG-algebra. To circumvent this we consider
the sheaf-theoretic pullback functor $\pi^{-1}$
(recall that $\pi^*G \cong \pi^{-1}G\otimes_{\pi^{-1}\CO_Y}\CO_\TY$ by definition),
and define
$$
\Psi:\D^-(Y,\CA) \to \D^-(\TY),\qquad
G \mapsto \pi^{-1}G \otimes_{\pi^{-1}(\CA)} E.
$$
Note that the functor $\pi^{-1}$ is exact, hence $\pi^{-1}\CA$ is a sheaf of algebras on $\TY$
and $\Psi$ is right-exact. Note also that the adjunction morphism
$\pi^{-1}\CA = \pi^{-1}\pi_*(\CEnd E) \to \CEnd E$ provide $E$ with a structure
of a left $\pi^{-1}\CA$-module (commuting with the right $\CO_\TY$-module structure),
giving a sense to the above tensor product.

\begin{lemma}\label{gus}
The functor $\Psi:\D^-(Y,\CA) \to \D^-(\TY)$ is left adjoint to
$\Phi:\D^-(\TY) \to \D^-(Y,\CA)$. Moreover, the composition
$\Phi\circ\Psi:\D^-(Y,\CA) \to \D^-(Y,\CA)$
is isomorphic to the identity functor $\Phi\circ\Psi \cong \id$.
\end{lemma}
\begin{proof}
The adjointness is straightforward
\begin{multline*}
\Hom(\Psi(G),F) =
\Hom_{\CO_\TY}(\pi^{-1}G \otimes_{\pi^{-1}(\CA)} E,F) \cong \\ \cong
\Hom_{\pi^{-1}(\CA)}(\pi^{-1}(G),F \otimes_{\CO_\TY} E^*) \cong 
\Hom_{\CA}(G,\pi_*(F \otimes_{\CO_\TY} E^*)) =
\Hom_{\CA}(G,\Phi(F)).
\end{multline*}
Similarly, we have
$$
\Phi\circ\Psi(G) =
\pi_*(\pi^{-1}G \otimes_{\pi^{-1}(\CA)} E \otimes_{\CO_\TY} E^*) \cong 
\pi_*(\pi^{-1}G \otimes_{\pi^{-1}(\CA)} \CEnd E) \cong
G \otimes_{\CA} \pi_*(\CEnd E) \cong
G \otimes_{\CA} \CA \cong
G
$$
which shows that the composition is the identity functor.
\end{proof}


\begin{lemma}
If $P \in \D^-(Y)$, $F \in \D^-(\TY)$, $G \in \D^-(Y,\CA)$ then we have
$\Phi(\pi^* P\otimes_{\CO_\TY} F) \cong P\otimes_{\CO_Y} \Phi(F)$ and
$\Psi(P \otimes_{\CO_Y} G) \cong \pi^*P \otimes_{\CO_\TY} \Psi(G)$.
\end{lemma}
\begin{proof}
The first isomorphism follows from the projection formula
$\pi_*(\pi^* P \otimes F \otimes E^*) \cong P\otimes \pi_*(F\otimes E^*)$.
The second is also quite easy
\begin{multline*}
\Psi(P\otimes_{\CO_\TY} G) =
\pi^{-1}(P\otimes_{\CO_\TY} G)\otimes_{\pi^{-1}\CA} E \cong
\pi^{-1}P\otimes_{\pi^{-1}\CO_Y} (\pi^{-1}G\otimes_{\pi^{-1}\CA} E) \cong
\\ \cong
\pi^{-1}P\otimes_{\pi^{-1}\CO_Y}\CO_\TY\otimes_{\CO_\TY}\Psi(G) \cong
\pi^*P\otimes_{\CO_Y}\Psi(G).
\end{multline*}
\end{proof}

We will also need the following property of the functor $\Psi$.

\begin{lemma}\label{psijs}
If $G \in \D^-(Z)$ and $j_*G \in \D^-(Y,\CA)$ then
$\Psi(j_*G) \cong i_*F$ for some $F \in \D^-(\TZ)$.
\end{lemma}
\begin{proof}
Let $G^\bullet$ be a bounded above complex of sheaves of $j^{-1}\CA$-modules on $Z$ representing $G$.
Choose any flat over $\CA$ resolution $\CE^\bullet$ of $E$ in the category of $(\pi^{-1}\CA\!-\!\CO_\TY)$-bimodules.
Then $\Psi(j_*G^\bullet)$ can be represented by the total complex of a bicomplex
$$
\CG^{\bullet,\bullet} = \pi^{-1}j_*G^\bullet\otimes_{\pi^{-1}\CA}\CE^\bullet.
$$
Let us show that each term of this bicomplex (which a priori is a sheaf of $\CO_\TY$-modules)
is supported (scheme-theoretically) on $\TZ$. Indeed, since $\pi^{-1}\CO_Y$ is a central subalgebra
in $\pi^{-1}\CA$, the $\pi^{-1}\CO_Y$-module structure on the term $\CG^{p,q}$ comes from that
on $\pi^{-1}j_*G^p$, hence $\CG^{p,q}$ is annihilated by $\pi^{-1}\CJ_Z$, where $\CJ_Z \subset \CO_Y$
is the sheaf of ideals of $Z$ in $Y$. Since by~\eqref{pisjz} the ideal generated by $\pi^{-1}\CJ_Z$ in $\CO_\TY$
is $\CJ_\TZ$, we conclude that $\CG^{p,q}$ is supported
scheme-theoretically on $\TZ$, in the other words $\CG^{p,q} = i_*\CF^{p,q}$,
where $\CF^{p,q}$ is a coherent sheaf of $\CO_\TZ$-modules.
Since the pushforward functor $i_*$ is fully faithful on the level of abelian categories,
the sheaves $\CF^{p,q}$ form a bicomplex on $\TZ$ and we have
$\CG^{\bullet,\bullet} \cong i_*\CF^{\bullet,\bullet}$,
hence $\Psi(j_*G) \cong i_*F$, where $F$ is an object of $\D^-(\TZ)$
represented by (the total complex of) a bicomplex $\CF^{\bullet,\bullet}$.
\end{proof}

Let $\CB_k^-(-k) \subset \D^-(\TZ)$ and $\TD^- \subset \D^-(\TY)$ be the
negative completions of $\CB_k \subset \D^b(\TZ)$ and $\TD \subset \D^b(\TY)$.
By lemma~\ref{dminsod} and proposition~\ref{caf} we have
the following semiorthogonal decompositions
$$
\begin{array}{l}
\D^-(\TZ) = \langle \CB_{m-1}^-(1-m),\CB_{m-2}^-(2-m),\dots,\CB_1^-(-1),\CB_0^-\rangle,\\
\D^-(\TY) = \langle i_*(\CB_{m-1}^-(1-m)),i_*(\CB_{m-2}^-(2-m)),\dots,i_*(\CB_1^-(-1)),\TD^-\rangle.
\end{array}
$$

\begin{lemma}\label{phi0}
Let $F \in \D^-(\TZ)$.
Then $\Phi(i_*F) = 0$ if and only if $F \in \langle \CB_{m-1}^-(1-m),\dots,\CB_1^-(-1) \rangle$.
\end{lemma}
\begin{proof}
We have
$$
\Phi(i_*F) =
\pi_*(i_*F \otimes E^*) =
\pi_*i_*(F\otimes i^*E^*) =
j_*p_*\RCHom(i^*E,F),
$$
so $\Phi(i_*F) = 0$ is equivalent to $p_*\RCHom(i^*E,F) = 0$.
Note that for any $G \in \D^-(Z)$ we have
$$
\Hom(G,p_*\RCHom(i^*E,F)) \cong
\Hom(p^*G,\RCHom(i^*E,F)) \cong
\Hom(p^*G\otimes i^*E,F),
$$
therefore $p_*\RCHom(i^*E,F) = 0$ is equivalent to $\Hom(i^*E\otimes p^*(\D^-(Z)),F) = 0$.
Further, \eqref{ldtz1} implies that $i^*E \otimes p^*(\D^-(Z)) \supset \CB_0^-$,
hence $\Hom(i^*E\otimes p^*(\D^-(Z)),F) = 0$ implies $F \in (\CB_0^-)^\perp$,
and it remains to note that $(\CB_0)^\perp = \langle \CB_{m-1}^-(1-m),\dots,\CB_1^-(-1) \rangle$.
\end{proof}

\begin{proposition}\label{minus}
The functor $\Psi$ gives an equivalence of categories
$$
\Psi: \D^-(Y,\CA) \overset{\sim}\longrightarrow \TD^-.
$$
\end{proposition}
\begin{proof}
We know by lemma~\ref{phi0} that $\Phi$ vanishes on
$\lan i_*(\CB_{m-1}^-(1-m)),i_*(\CB_{m-2}^-(2-m)),\dots,i_*(\CB_1^-(-1))\ran$.
It follows by adjunction that the functor $\Psi$ takes $\D^-(Y,\CA)$ to $\TD^-$.
Moreover, by lemma~\ref{gus} it is fully faithful.
Let us check that it gives an equivalence $\D^-(Y,\CA) \cong \TD^-$.
In other words, let us check that we have the following semiorthogonal decomposition
\begin{equation}\label{dbmty}
\D^-(\TY) = \langle i_*(\CB_{m-1}^-(1-m)),i_*(\CB_{m-2}^-(2-m)),\dots,i_*(\CB_1^-(-1)),\D^-(Y,\CA)\rangle.
\end{equation}

First of all, we will show that the sheaf $i_*i^*E$ is contained in the RHS of~\eqref{dbmty}.
For this we note that
$$
\Phi(i_*i^*E) = \pi_*(i_*i^*E \otimes E^*) \cong \pi_*i_*i^*\CEnd E \cong j_*p_*i^*\CEnd E
$$
and by lemma~\ref{psijs} we deduce that $\Psi(\Phi(i_*i^*E))$ lies
in the subcategory $i_*(\D^-(\TZ))$ of $\D^-(\TY)$.

%

Now, consider the adjunction morphism
$\Psi\Phi(i_*i^*E) \to i_*i^*E$ and its cone.
Since this map evidently comes from a map in $\D^-(\TZ)$
it follows that there exists $F\in\D^-(\TZ)$ such that
we have the following distinguished triangle
\begin{equation}\label{triang}
\Psi\Phi(i_*i^*E) \to i_*i^*E \to i_*F.
\end{equation}
Applying the functor $\Phi$ and taking into account lemma~\ref{gus}
we deduce that $\Phi(i_*F) = 0$. By lemma~\ref{phi0} it follows that
$F \in \langle\CB^-_{m-1}(1-m),\dots,\CB^-_1(-1)\rangle$, hence
$i_*F \in \langle i_*\CB^-_{m-1}(1-m),\dots,i_*\CB^-_1(-1)\rangle$.
Thus we see that the third term of~\eqref{triang} is contained in the RHS of~\eqref{dbmty}.
On the other hand, the first term of~\eqref{triang} is contained in the RHS of~\eqref{dbmty}
by definition. Therefore the second term, $i_*i^*E$, is also contained in the RHS of~\eqref{dbmty}.

Further, take any $P \in \D^-(Y)$ and note that
$\pi^*P\otimes\Psi(\Phi(i_*i^*E)) \cong \Psi(P\otimes\Phi(i_*i^*E))$
is contained in $\Psi(\D^-(Y,\CA))$. On the other hand,
$\pi^*P \otimes i_*F \cong i_*(i^*\pi^*P \otimes F)$
and $\Phi(\pi^*P \otimes i_*F) \cong P\otimes \Phi(i_*F) = 0$,
hence $\pi^*P \otimes i_*F \in \langle i_*(\CB^-_{m-1}(1-m)),i_*(\CB^-_{m-2}(2-m)),\dots,i_*(\CB^-_1(-1))\rangle$.
Tensoring~\eqref{triang} by $\pi^*P$ we conclude that
$\pi^* P \otimes i_*i^*E$ is also contained in the RHS of~\eqref{dbmty}
for any $P\in\D^-(Y)$.

Now let us check that the RHS of~\eqref{dbmty} coincides with the LHS.
Since the RHS is evidently left admissible, it suffices to check
that the right orthogonal to the RHS is zero.
Assume that $G\in\D^-(\TY)$ lies in the right orthogonal to the RHS.
Since $\pi^*P\otimes i_*i^*E$ is contained in the RHS for any $P\in\D^-(Y)$
we deduce that
$$
0 =
\Hom(\pi^*P\otimes i_*i^*E,G) =
\Hom(\pi^*P,\RCHom(i_*i^*E,G)) =
\Hom(P,\pi_*\RCHom(i_*i^*E,G)),
$$
hence $\pi_*\RCHom(i_*i^*E,G) = 0$ by lemma~\ref{perforth}.
But
$$
0 =
\pi_*\RCHom(i_*i^*E,G) \cong
\pi_*i_*\RCHom(i^*E,i^!G) \cong
\pi_*\RCHom(E,i_*i^!G) \cong
\Phi(i_*i^!G)
$$
hence $i^!G \in \langle\CB_{m-1}^-(1-m),\dots,\CB_1^-(-1)\rangle$ by lemma~\ref{phi0}.
On the other hand, for all $1\le k\le m-1$ we have
$0 = \Hom(i_*\CB_k^-(-k),G) = \Hom(\CB_k^-(-k),i^!G)$,
hence $i^!G \in \langle\CB_{m-1}^-(1-m),\dots,\CB_1^-(-1)\rangle^\perp$.
Combining, we deduce that $i^!G = 0$.

Since $i:\TZ \to \TY$ is a closed embedding, $i^!G = 0$ implies that $G$ is supported
on the complement of $\TZ$ in~$\TY$. Finally we note that $\pi$ is identity on this complement,
hence
$$
0 = \Hom(\Psi(\D^-(Y,\CA)),G) \cong
\Hom(\D^-(Y,\CA),\Phi(G)) \cong
\Hom(\D^-(Y,\CA),\pi_*(G\otimes E^*))
$$
implies $0 = \pi_*(G\otimes E^*) = G\otimes E^*$
hence $G = 0$ since $E$ is a vector bundle.
\end{proof}

Thus we see that the functor $\Psi: \D^-(Y,\CA) \to \TD^-$ is an equivalence.
Our further goal is to check that this equivalence restricts to an equivalence
of $\D^b(Y,\CA)$ and $\TD$. As an intermediate step we check the following

\begin{lemma}
The functor $\Phi$ takes $\D^b(\TY)$ to $\D^\perf(Y,\CA)$.
\end{lemma}
\begin{proof}
Let $F \in \D^b(\TY)$. It is clear that $\Phi(F)$ is bounded.
So, by proposition~\ref{perfect} it suffices to check that
$\Phi(F)$ has finite $\Ext$-amplitude. Indeed, take any $G \in \D^{\le n}(Y,\CA)$.
Then it is clear that $\Psi(G) \in \D^{\le n}(\TY)$.
On the other hand, the functor $\Psi\Phi:\D^-(\TY) \to \D^-(\TY)$
is the projection to the $\TD^-$ component of the decomposition~\eqref{dbmty}.
Since $F$ is bounded lemma~\ref{dbdm} shows that $\Psi\Phi(F)$ is also bounded.
Since $\TY$ is smooth it follows from lemma~\ref{hub_loc} that it has
finite $\Ext$-amplitude $a = a(\Psi\Phi(F))$. We conclude by
$$
\RHom(\Phi(F),G) = \RHom(\Psi\Phi(F),\Psi(G)) \in \D^{\le n+a}(\kk),
$$
hence $\Phi(F)$ has finite $\Ext$-amplitude.
\end{proof}

We conclude the proof of theorem~\ref{th2} by the following

\begin{proposition}
The functor $\Psi$ gives an equivalence $\D^b(Y,\CA) \cong \TD$.
In particular, the sheaf of algebras $\CA$ has finite homological dimension.
\end{proposition}
\begin{proof}
Assume that $G\in\D^b(Y,\CA)$ and $\Psi(G)$ is unbounded from below.
Then by lemma~\ref{dminus} there exists $F\in\D^b(\TY)$ such that
$\RHom(F,\Psi(G))$ is unbounded from below. Applying to $F$
the left adjoint functor to the embedding $\TD \to \D^b(\TY)$
(recall that $\TD$ is left admissible, see remark~\ref{tdadm})
we obtain an object $F'\in\TD$ such that $\RHom(F',\Psi(G))$ is unbounded from below.
Since $\Phi$ is an equivalence,
$$
\RHom(F',\Psi(G)) \cong
\RHom(\Phi(F'),\Phi(\Psi(G))) =
\RHom(\Phi(F'),G),
$$
and we deduce that $\RHom(\Phi(F'),G)$ is unbounded from below.
But as we already proved $\Phi$ takes $\TD$ to $\D^\perf(Y,\CA)$,
hence $\Phi(F') \in \D^\perf(Y,\CA)$. Thus we got a contradiction:
$\RHom$ between a perfect complex and a bounded object is always bounded.

We have checked that $\Psi$ takes $\D^b(Y,\CA)$ to $\TD$.
On the other hand, $\Phi$ takes $\TD$ to $\D^b(Y,\CA)$.
So, they give equivalences $\D^b(Y,\CA) \cong \TD$.
Since $\TD$ is a semiorthogonal component of $\D^b(\TY)$
and $\TY$ is smooth we conclude that $\CA$ has finite homological dimension.
\end{proof}

\begin{remark}
The functor $\D^b(Y,\CA) \to \D^b(Y)$ providing $\D^b(Y,\CA)$ with a structure of
a categorical resolution of singularities of $\D^b(Y)$ is {\em not}\/ the forgetting
of the $\CA$-module structure functor. Actually, it is given by $F \mapsto F\otimes_\CA\pi_*E$
(and its left adjoint functor $\D^\perf(Y) \to \D^b(Y,\CA)$ is given by
$G \mapsto G\otimes_{\CO_Y}\pi_*E^*$).
\end{remark}

\section{Functoriality}\label{funct}

Assume that $Y$, $Z$, $\TY$, and $\TZ$ are as in section~\ref{ncr}, that is
$Y$ is a rational singularity, $\pi:\TY \to Y$ is a resolution of singularities
such that its exceptional locus $\TZ$ is an irreducible divisor, $Z = \pi(\TZ)$;
$i:\TZ \to \TY$ denotes the embedding, and $p:\TZ \to Z$ denotes the projection.

Let $\phi:Y' \to Y$ be a change of base.
Put $Z' = Z\times_Y Y'$, $\TY' = \TY\times_Y Y'$
and $\TZ' = \TZ\times_Y Y'$. Denote $\pi':\TY' \to Y'$, $p':\TZ' \to Z'$, and $i':\TZ'\to\TY'$.

We restrict in this section to the following class of base changes. We assume that
\begin{itemize}
\item $\phi$ can be factored as a locally complete intersection closed embedding
followed by a smooth projective morphism (in particular, $\phi$ has finite $\Tor$-dimension);
\item both $Y'$ and $Y$ are Cohen--Macaulay and both $Z'$, $\TZ'$ and $\TY'$ have expected dimension.
\end{itemize}
For $\phi$ satisfying these properties it is easy to show that $\pi'$ is birational,
$\TZ'$ is its exceptional divisor and $Z' = \pi'(\TZ')$. Moreover, the following base-change isomorphisms
follow from~\cite{K1}, section~2.4
\begin{equation}\label{bcipi}
\hspace{-10pt}
\vcenter{\xymatrix{
\TZ' \ar[r]^\phi \ar[d]_{i'} & \TZ \ar[d]_i \\
\TY' \ar[r]^\phi \ar[d]_{\pi'} & \TY \ar[d]_\pi \\
Y' \ar[r]^\phi & Y
}}
\quad
\vcenter{\xymatrix{
\TZ' \ar[r]^\phi \ar[d]_{p'} & \TZ \ar[d]_p \\
Z' \ar[r]^\phi \ar[d]_{j'} & Z \ar[d]_j \\
Y' \ar[r]^\phi & Y
}}
\quad
\begin{array}{ll}
i'_*\phi^* \cong \phi^* i_*:\D^b(\TZ) \to \D^b(\TY'), &
\phi_*{i'}^* \cong i^*\phi_*:\D^b(\TY') \to \D^b(\TZ), \\[7pt]
p'_*\phi^* \cong \phi^* p_*:\D^b(\TZ) \to \D^b(Z'), &
\phi_*{p'}^* \cong p^*\phi_*:\D^b(Z') \to \D^b(\TZ), \\[12pt]
\pi'_*\phi^* \cong \phi^*\pi_*:\D^b(\TY) \to \D^b(Y'), &
\phi_*{\pi'}^* \cong \pi^*\phi_*:\D^b(Y') \to \D^b(\TY),\\[7pt]
\j'_*\phi^* \cong \phi^*j_*:\D^b(Z) \to \D^b(Y'), &
\phi_*{j'}^* \cong j^*\phi_*:\D^b(Y') \to \D^b(Z).
\end{array}
\end{equation}

Assume that we are given a dual Lefschetz decomposition~\eqref{ldtz}.
Let $\CB'_k = \langle \phi^*(\CB_k)\otimes {p'}^*(\D^\perf(Z')) \rangle \subset \D^b(\TZ')$.
Then the categories $\CB'_{m-1}(1-m),\CB'_{m-2}(2-m),\dots,\CB'_1(-1),\CB'_0$ form a semiorthogonal
collection in $\D^b(\TZ')$. We assume that this collection generates $\D^b(\TZ')$, i.e.
$$
\D^b(\TZ') = \langle \CB'_{m-1}(1-m),\CB'_{m-2}(2-m),\dots,\CB'_1(-1),\CB'_0 \rangle
$$
is a semiorthogonal decomposition.
\begin{remark}
This assumption automatically holds if each of the categories $\CB_k$
is generated (over $Z$) by an exceptional over $Z$ collection of vector bundles on $\TZ$.
\end{remark}

Assume also that $\CB_0$ is $Z$-linear. Then
$$
\phi_*(\phi^*(\CB_0)\otimes {p'}^*(\D^\perf(Z'))) \subset
\CB_0\otimes \phi_*{p'}^*(\D^\perf(Z'))) \subset
\CB_0\otimes p^*\phi_*(\D^\perf(Z'))) \subset
\CB_0\otimes p^*(\D^\perf(Z))) \subset
\CB_0
$$
(we have used the base-change~\eqref{bcipi} and finiteness of $\Tor$-dimension of $\phi$), hence
$\phi_*(\CB'_0) \subset \CB_0$.

Denote by $\TD' \subset \D^b(\TY')$ the full subcategory of $\D^b(\TY')$
consisting of all $F \in \D^b(\TY')$ such that $i^*F \in \CB'_0$.
If $\TY'$ is smooth then $\TD'$ is a categorical resolution of singularities
of $\D^b(Y')$ by theorem~\ref{th1}. However, in general $\TY'$ can be singular.
Anyway, we have the following

\begin{proposition}
We have a semiorthogonal decomposition
$$
\D^b(\TY') = \lan i_*(\CB'_{m-1}(1-m)),i_*(\CB'_{m-2}(2-m)),\dots,i_*(\CB'_1(-1)),\TD' \ran.
$$
Moreover, the functors $\phi^*:\D^b(\TY) \to \D^b(\TY')$ and $\phi_*:\D^b(\TY') \to \D^b(\TY)$
take $\TD$ to $\TD'$ and $\TD'$ to $\TD$ respectively, and we have the following base change isomorphisms
$$
\pi'_*\phi^* \cong \phi^*\pi_*:\TD \to \D^b(Y'),\qquad
\phi_*{\pi'}^* \cong \pi^*\phi_*:\D^b(Y') \to \TD.
$$
\end{proposition}
\begin{proof}
The semiorthogonal decomposition follows from proposition~\ref{isld}.
Further, if $F \in \TD$ then we have
${i'}^*\phi^*F \cong \phi^*i^*F \in \phi^*(\CB_0) \subset \CB'_0$,
hence $\phi^*F \in \TD'$. Similarly, using the base-change~\eqref{bcipi} we deduce
$i^*\phi_*F \cong \phi_*{i'}^*F \in \phi_*(\CB'_0) \subset \CB_0$.
\end{proof}

Further, assume that there is a vector bundle $E$ on $\TY$ such that~\eqref{ldtz1},
\eqref{defa} and~\eqref{pisjz} holds. Denote
$$
E' = \phi^*E,
$$
a vector bundle on $Y'$.

\begin{lemma}
The vector bundle $E'$ on $Y'$ enjoys the corresponding
properties~\eqref{ldtz1}, \eqref{defa} and~\eqref{pisjz}.
\end{lemma}
\begin{proof}
First of all, \eqref{ldtz1} holds by definition of $\CB'_0$ and $E'$.
Further, it is clear that $\CA' = \pi'_*\CEnd E' \in \D^{\ge 0}(Y)$.
On the other hand, $\CA' = \pi'_*\CEnd E' \cong \pi'_*\phi^*\CEnd E \cong
\phi^*\pi_*\CEnd E = \phi^*\CA \in \D^{\le 0}(Y)$. Hence $\CA' \in \Coh(Y)$
and~\eqref{defa} holds.
Finally, we have $\phi^*j_*\CO_Z = j'_*\phi^*\CO_Z = j'_*\CO_{Z'}$, hence
$\phi^*\CJ_Z \cong \CJ_{Z'}$ and similarly $\phi^*\CJ_\TZ = \CJ_{\TZ'}$.
It follows that ${\pi'}^*\CJ_{Z'} = {\pi'}^*\phi^*\CJ_Z = \phi^*\pi^*\CJ_Z$,
and since ${\pi'}^{-1}\CJ_{Z'}\cdot\CO_{\TY'}$ is the image of ${\pi'}^*\CJ_{Z'}$ in $\CO_{\TY'} = \phi^*\CO_\TY$,
it follows that ${\pi'}^{-1}\CJ_{Z'}\cdot\CO_{\TY'} = \phi^*(\pi^{-1}\CJ_Z\cdot\CO_\TY) = \phi^*\CJ_\TZ = \CJ_{\TZ'}$.
\end{proof}

Define the functor $\Phi':\D^b(\TY') \to \D^b(Y',\CA')$ as $\Phi'(F) = \pi'_*(F\otimes {E'}^*)$.

\begin{theorem}
The functor $\Phi'$ induces an equivalence $\TD' \cong \D^b(Y',\CA')$ and we have the following base-change isomorphisms
$$
\Phi'\phi^* \cong \phi^*\Phi:\TD \to \D^b(Y',\CA'),\qquad
\Phi'\phi^! \cong \phi^!\Phi:\TD \to \D^b(Y',\CA'),\qquad
\phi_*\Phi' \cong \Phi\phi_*:\TD' \to \D^b(Y,\CA).
$$
\end{theorem}
\begin{proof}
First, let us check the base change
$$
\Phi'(\phi^*F) =
\pi'_*(\phi^* F\otimes {E'}^*) \cong
\pi'_*(\phi^* F\otimes \phi^* E^*) \cong
\pi'_*\phi^* (F\otimes E^*) \cong
\phi^*\pi_* (F\otimes E^*) =
\phi^*\Phi(F).
$$
and similarly for $\phi^!$ instead of $\phi^*$.
Further,
$$
\phi_*\Phi'(F) =
\phi_*\pi'_*(F\otimes {E'}^*) \cong
\phi_*\pi'_*(F\otimes \phi^* E^*) \cong
\pi_*\phi_*(F\otimes \phi^* E^*) \cong
\pi_*(\phi_*F\otimes E^*) =
\Phi(\phi_* F).
$$
It remains to check the equivalence $\TD' \cong \D^b(Y',\CA')$.
Recall that we assumed that $\phi$ can be factored as a locally complete intersection closed embedding
followed by a smooth projective map. By functoriality, it suffices to prove the equivalence for locally
complete intersection closed embeddings and for smooth projective morphisms.

If $\phi$ is smooth and projective, then in particular $\TY'$ is smooth, hence we can apply theorem~\ref{th2}.
If $\phi$ is a closed embedding, consider the functor $\Psi':\D^-(Y',\CA') \to \D^-(\TY')$ defined by
$\Psi'(G) = {\pi'}^{-1}G\otimes_{{\pi'}^{-1}\CA'}E'$. The same arguments as in the proof of lemma~\ref{gus}
and proposition~\ref{minus} show that $\Psi'$ is left adjoint to $\Phi'$ and that they give
mutually inverse equivalences ${\TD'}{}^-$ and $\D^-(Y',\CA')$. Therefore, it suffices to check
that $\Psi'$ preserves boundedness. For this we note that $\phi_*\Psi' \cong \Psi\phi_*$
(this follows from the isomorphism $\Phi'\phi^! \cong \phi^!\Phi$ by taking left adjoints of both sides).
Therefore for any $G\in\D^b(Y',\CA')$ we have $\phi_*\Psi'(G) \cong \Psi\phi_*(G)$ is bounded,
hence $\Psi'(G)$ is also bounded, since $\phi$ is a closed embedding.
\end{proof}

\section{Cones}\label{CONES}

Let $Y$ be an affine cone over a smooth algebraic variety $X$ with respect to a line bundle $L$.
Then the singular locus $Z$ of $Y$ is the vertex of the cone, $Z = \Spec\kk$.
The blowup $\pi:\TY \to Y$ of $Y$ is smooth. Moreover, $\TY$ is isomorphic
to the total space of the line bundle $L^{-1}$ on $Y$,
and the exceptional divisor $\TZ$ of the blowup coincides
with the zero section of $L^{-1}$.
Thus
$$
\TZ \cong X,
\qquad\text{ }\qquad
\CN^*_{\TZ/\TY} \cong L
$$
and, moreover, the condition~\eqref{pisjz} is satisfied.

As we have seen any (dual) Lefschetz decomposition of $\D^b(\TZ) = \D^b(X)$
with respect to $L$ gives rise to a categorical resolution of $\D^b(Y)$.
Furthermore, this resolution is a noncommutative resolution of singularities
if the corresponding Lefschetz decomposition is generated by an appropriate
vector bundle $E$ on $\TY$ which is tilting over $Y$.
Now we will give a useful criterion of tiltingness.

Let $q:\TY \to X$ be the natural projection
(recall that $\TY$ is isomorphic to the total space of the line bundle $L^{-1}$ on $X$).
\begin{proposition}\label{criterion}
Let $F$ be a vector bundle on $X$.
The vector bundle $E = q^*F$ on $\TY$ is tilting over $Y$
if and only if we have
\begin{equation}\label{tilting}
H^{>0}(X,\CEnd F\otimes L^t) = 0
\qquad\text{for all $t\ge 0$}.
\end{equation}
\end{proposition}
\begin{proof}
Since $Y$ is affine $R^{>0}\pi_*(\CEnd E) = 0$ is equivalent to
$H^{>0}(Y,\pi_*(\CEnd E)) = H^{>0}(\TY,\CEnd E) = 0$.
On the other hand,
\begin{multline*}
H^\bullet(\TY,\CEnd E) \cong
H^\bullet(X,q_*\CEnd E) \cong
H^\bullet(X,q_*q^*\CEnd F) \cong
\\ \cong
H^\bullet(X,\CEnd F \otimes q_*\CO_\TY) \cong
H^\bullet(X,\CEnd F \otimes (\oplus_{t=0}^\infty L^t)) \cong
\bigoplus_{t=0}^\infty H^\bullet(X,\CEnd F \otimes L^t).
\end{multline*}
and the proposition follows.
\end{proof}

\subsection{Veronese cones}

Let $Y$ be a cone over $\PP^n$ in the degree $d$ Veronese embedding with $d \le n+1$
(in other words, $Y$ is the weighted projective space $\PP(1^{n+1},d)$), and $Z \subset Y$
be the singular point.
The exceptional divisor $\TZ$ of the blowup $\pi:\TY \to Y$ of $Z$ coincides with $\PP^n$
and the conormal bundle $\CN^*_{\TZ/\TY}$ coincides with $\CO_\TZ(d)$.
Define integers $m\ge 1$ and $1\le r\le d$ from the equality $n+1 = (m-1)d + r$.
Consider the following dual Lefschetz decomposition of $\D^b(\TZ)$
$$
\D^b(\TZ) = \lan \CB_{m-1}(-d(m-1)),\dots,\CB_1(-d),\CB_0\ran,
$$
where $\CB_0 = \CB_1 = \dots = \CB_{m-2} = \lan \CO_\TZ(1-d),\dots,\CO_\TZ(-1),\CO_\TZ\ran$,
and $\CB_{m-1} = \lan \CO_\TZ(1-r),\dots,\CO_\TZ(-1),\CO_\TZ\ran$.
We conclude that by theorem~\ref{th1} the category
$$
\TD = {}^\perp(i_*\CO_\TZ(-n),\dots,i_*\CO_\TZ(-d)) \subset \D^b(\TY)
$$
is a categorical resolution of singularities.
It is crepant if and only if $r = d$, that is, if $d$ divides $n+1$.

Moreover, taking $E = q^*\CO_\TZ \oplus \dots \oplus q^*\CO_\TZ(d-1)$
we deduce that~\eqref{ldtz1} holds. It is also easy to see that
the vector bundle $E$ on $\TY$ is tilting over $Y$ by proposition~\ref{criterion},
hence by theorem~\ref{th2}
$$
\TD = \D^b(Y,\CA)
$$
is a noncommutative resolution with $\CA = \pi_*\CEnd E$.
It is crepant if and only if $d$ divides $n+1$.

\subsection{Segre cones}

Let $Y$ be a cone over $\PP^{m-1}\times\PP^{m-1}$ in the Segre embedding, and $Z \subset Y$
be the singular point.
The exceptional divisor $\TZ$ of the blowup $\pi:\TY \to Y$ of $Z$ coincides with $\PP^{m-1}\times\PP^{m-1}$
and the conormal bundle $\CN^*_{\TZ/\TY}$ coincides with $\CO_\TZ(1,1)$.
Consider the following dual Lefschetz decomposition of $\D^b(\TZ)$
$$
\D^b(\TZ) = \lan \CB_{m-1}(1-m,1-m),\dots,\CB_1(-1,-1),\CB_0\ran,
$$
where $\CB_0 = \CB_1 = \dots = \CB_{m-1} = \lan \CO_\TZ(0,1-m),\dots,\CO_\TZ(0,-1),\CO_\TZ\ran$.
We conclude that by theorem~\ref{th1} the category
$$
\TD = {}^\perp(i_*\CB_0(1-m,1-m),\dots,i_*\CB_0(-1,-1)) \subset \D^b(\TY)
$$
is a categorical crepant resolution of singularities.

Moreover, if $E = q^*\CO_\TZ \oplus q^*\CO_\TZ(0,-1) \oplus \dots \oplus q^*\CO_\TZ(0,1-m)$
then $\CB_0 = \lan i^*E \ran$ and it is easy to see that $E$ is tilting over $Y$ by proposition~\ref{criterion},
hence by theorem~\ref{th2}
$$
\TD = \D^b(Y,\CA)
$$
is a noncommutative resolution with $\CA = \pi_*\CEnd E$.

\subsection{Anticanonical cones}

Let $Y$ be the cone over a Fano manifold in the anticanonical embedding, and $Z \subset Y$ be the singular point.
The exceptional divisor $\TZ$ of the blowup $\pi:\TY \to Y$ of $Z$ is then
the initial Fano manifold, and the conormal bundle $\CN^*_{\TZ/\TY}$
coincides with $\omega^{-1}_\TZ$.
Then we can consider the stupid Lefschetz decomposition $\CB_0 = \D^b(\TZ)$,
and see that the category $\D^b(\TY)$ is a categorical crepant resolution of $\D^b(Y)$
(actually, $\TY$ is a crepant resolution of $Y$ in this case).
Moreover, assume that $\TZ$ admits a strong exceptional collection of vector bundles $\{E_s\}$.
Then we take $E = \oplus q^*E_s$. It is clear that $\CB_0 = \D^b(\TZ) = \lan i^*E \ran$
and $E$ is tilting over $Y$ by proposition~\ref{criterion}.
Therefore by theorem~\ref{th2}
$$
\TD = \D^b(Y,\CA)
$$
is a noncommutative crepant resolution with $\CA = \pi_*\CEnd E$.

\subsection{Cones over Grassmannians}\label{grass}

Let $Y$ be the cone over a Grassmannian of lines $\Gr(2,m)$ in the Pl\"ucker embedding,
and $Z \subset Y$ is the singular point.
The exceptional divisor $\TZ$ of the blowup $\pi:\TY \to Y$ of $Z$ then
coincides with the Grassmannian $\Gr(2,m)$ and the conormal bundle $\CN^*_{\TZ/\TY}$
coincides with $\CO_\TZ(1)$. Let $\CU$ denote the tautological rank $2$ bundle on $\TZ = \Gr(2,m)$.
Let $k = \lfloor \frac{m-1}2 \rfloor$.
Then we can consider the following Lefschetz decomposition
$$
\begin{array}{c}
\D^b(\TZ) = \lan \CB_{m-1}(1-m),\dots,\CB_1(-1),\CB_0 \ran, \qquad \text{where}\\
\begin{array}{ll}
\CB_0 = \CB_1 = \dots = \CB_{m-1} = \langle\CO_\TZ,\CU^*,\dots,S^k\CU^*\rangle \quad &
\text{if $m = 2k + 1$ is odd},\smallskip\\
\left\{
\begin{array}{l}
\CB_0 = \CB_1 = \dots = \CB_{k} = \langle\CO_\TZ,\CU^*,\dots,S^k\CU^*\rangle,\\
\CB_{k+1} = \CB_{k+2} = \dots = \CB_{2k+1} = \langle\CO_\TZ,\CU^*,\dots,S^{k-1}\CU^*\rangle,
\end{array}
\right.
\quad &
\text{if $m = 2k + 2$ is even}.
\end{array}
\end{array}
$$
We conclude that the category
$$
\TD = {}^\perp(i_*\CB_{m-1}(1-m),\dots,i_*\CB_1(-1)) \subset \D^b(\TY)
$$
is a categorical resolution of singularities. It is crepant iff $m$ is odd.

Moreover, if $E = q^*\CO_\TZ \oplus q^*\CU^* \oplus \dots \oplus q^*S^k\CU^*$
then $\CB_0 = \lan i^*E \ran$ and it is easy to see that $E$ is tilting over $Y$ by proposition~\ref{criterion}, hence
$$
\TD = \D^b(Y,\CA)
$$
is a noncommutative resolution with $\CA = \pi_*\CEnd E$. It is crepant iff $m$ is odd.

\section{The Pfaffian varieties}\label{PF}

Let $W$ be a vector space over $\kk$, $\dim W = n$.
Consider the projective space $\PP(\Lambda^2W^*)$ of skew-forms on $W$.
For each $0 \le t \le \lfloor n/2\rfloor$ consider the following
closed subset of $\PP(\Lambda^2W^*)$
$$
\Pf(2t,n) = \Pf(2t,W^*) = \PP(\{ \omega \in \Lambda^2W^*\ |\ \rank(\omega) \le 2t \}),
$$
where $\rank(\omega)$ is the rank of $\omega$ (the dimension of the image
of the map $W \to W^*$ induced by $\omega$).
We call $\Pf(2\lfloor n/2\rfloor-2,W^*)$ the {\sf Pfaffian variety},
and other $\Pf(2t,W^*)$ are called the {\sf generalized Pfaffian varieties}.
The ideal of the generalized Pfaffian variety $\Pf(2t,W^*)$ is generated
by the Pfaffians of all diagonal $(2t+2)\times(2t+2)$-minors of a skew-form,
hence the name.

\begin{example}
\begin{itemize}
\item $\Pf(0,W^*) = \emptyset$;
\item $\Pf(2,W^*) = \Gr(2,W^*)$;
\item $\Pf(2\lfloor n/2\rfloor,W^*) = \PP(\Lambda^2W^*)$;
\item for $t \ne 0,1,\lfloor n/2\rfloor$ the Pfaffian variety $\Pf(2t,W^*)$ is singular, the singular locus is the previous Pfaffian variety
$\Pf(2t-2,W^*) \subset \Pf(2t,W^*)$.
\end{itemize}
\end{example}

In this section we describe a noncommutative resolution of singularities
of the generalized Pfaffian variety $\Pf(4,W^*)$ for $n = \dim W \ge 6$.
So, put
$$
Y = \Pf(4,W^*),
\qquad
Z = \mathop{{\sf Sing}}(Y) = \Pf(2,W^*) = \Gr(2,W^*) = \Gr(n-2,W).
$$
Note that all skew-forms in $Y \setminus Z$ are of rank $4$, hence their kernels are $(n - 4)$-dimensional.
Similarly, all skew-forms in $Z$ are of rank $2$, and their kernels are $(n - 2)$-dimensional.
Let $\TY$ be the set of all pairs $(\omega,K)$, where $K$ is an $(n - 4)$-dimensional subspace
in $W$ and $\omega$ is a skew-form containing $K$ in the kernel.
More precisely,
$$
\TY = \PP_{\Gr(n - 4,W)}(\Lambda^2\CK^\perp),
$$
where $\CK \subset W\otimes\CO_{\Gr(n - 4,W)}$ is the tautological subbundle of rank $n - 4$,
and $\CK^\perp \subset W^*\otimes\CO_{\Gr(n - 4,W)}$ is the orthogonal to $\CK$ subbundle of rank~$4$.
Note that $\TY$ is irreducible and smooth.

Let us denote the pullbacks of the vector bundles $\CK$ and $\CK^\perp$
from $\Gr(n-4,W)$ to $\TY$ by the same letters. Note that the embedding
of vector bundles $\Lambda^2\CK^\perp \to \Lambda^2W^*\otimes\CO_\TY$
gives a map $\TY \to \PP(\Lambda^2W^*)$ which is birational onto $Y$.
Thus $\TY$ is a resolution of singularities of $Y$.

\begin{lemma}
The blowup of $Y$ in $Z$ coincides with $\TY$.
The exceptional divisor $\TZ$ of the blowup
coincides with $\Gr_{\Gr(n-4,W)}(2,\CK^\perp) \cong \Fl(n-4,n-2;W) \cong \Fl(2,4;W^*)$.
The projection $\TZ \to Z$ coincides with the natural projection of the flag variety
$\Fl(2,4;W^*)$ to $\Gr(2,W^*)$. The conormal bundle $\CN^*_{\TZ/\TY}$
is isomorphic to the pullback of $\CO_{\Gr(n-4,W)}(1)$.
\end{lemma}
\begin{proof}
Denote temporarily the blowup of $Y$ in $Z$ by $\HY$ and the exceptional divisor by $\HZ$.

Since $Z = \Gr(2,W^*)$ the sheaf of ideals $J_Z$ is generated by the space of quadrics in $\PP(\Lambda^2W^*)$
passing through $Z$, which coincides with $\Lambda^{n-4}W^* \cong \Lambda^4W \subset S^2(\Lambda^2 W)$,
therefore we have a map from the blowup $\HY = \Proj_{Y}(\oplus J_{Z}^n)$ to $\PP(\Lambda^{n-4}W)$.
It is clear that under this map any point on $Y\setminus Z$, which is a skew-form $\omega\in\Lambda^2W^*$
of rank $4$, goes to a bivector in $\Lambda^{n-4}W$ corresponding to the skew-form $\omega\wedge\omega \in \Lambda^4W^*$
under the isomorphism $\Lambda^4W^* \cong \Lambda^{n-4}W$. But for $\omega$ of rank $4$ this bivector is the decomposable
bivector corresponding to the $(n-4)$-dimensional subspace $\Ker\omega \subset W$. Thus the restriction of the map
$\HY \to \PP(\Lambda^{n-4}W)$ to $\HY\setminus\HZ = Y\setminus Z$ factors through $\Gr(n-4,W)\subset\PP(\Lambda^{n-4}W)$
hence the whole map factors as $\HY \to \Gr(n-4,W) \subset \PP(\Lambda^{n-4}W)$. The projections of the fibers of this
map to $Y$ can be identified with the spaces of all skew-forms containing given $(n-4)$-dimensional subspace $K\subset W$
in their kernels, thus we obtain a map $\HY \to \PP_{\Gr(n-4,W)}(\Lambda^2\CK^\perp) = \TY$.

Conversely, consider the map $\TY \to Y$. It is easy to see that the divisor of the pullback
of any quadric in $\PP(\Lambda^2W^*)$
passing through $Z$ is the union of the relative Grassmannian $\TZ = \Gr_{\Gr(n-4,W)}(\CK^\perp)$
and of the preimage of a hyperplane section of $\Gr(n-4,W)$, hence the pullback of the ideal $J_Z$,
generated by this quadrics is the ideal of the relative Grassmannian $\TZ$. It follows that the map
$\TY = \PP_{\Gr(n-4,W)}(\Lambda^2\CK^\perp) \to Y$ factors as $\TY \to \HY \to Y$.
The constructed maps $\HY \to \TY$ and $\TY \to \HY$ are mutually inverse, so the first claim of the lemma is proved.
Moreover, one can see that the other claims follow from the same arguments.
\end{proof}

Since $\TY$ is the blowup of $Y$ in $Z$, the condition~\eqref{pisjz} is satisfied.

Morphism $p:\TZ \to Z$ is smooth and its fibers are Grassmannians $\Gr(n-4,n-2) \cong \Gr(2,n-2)$.
We take the Lefschetz decomposition of $\D^b(\Gr(2,n-2))$ described in~\ref{grass}
and consider its relative version.

\begin{lemma}
Let
$$
\begin{array}{l l}
\left\{\arraycolsep=0.5pt\begin{array}{lllllllll}
\CB_0 &=& \CB_1 &=& \dots &=& \CB_{p-1} &=& \lan \CO_\TZ,\CK^*_{|\TZ},S^2\CK^*_{|\TZ},\dots,S^{p-1}\CK^*_{|\TZ} \ran_{\D^b(Z)},\\[5pt]
\CB_p &=& \CB_{p+1} &=& \dots &=& \CB_{2p-1} &=& \lan \CO_\TZ,\CK^*_{|\TZ},S^2\CK^*_{|\TZ},\dots,S^{p-2}\CK^*_{|\TZ} \ran_{\D^b(Z)}
\end{array}\right.
& \text{for $n = 2p+2$,}\\[10pt]
\hspace{10pt}\CB_0 = \CB_1 = \dots = \CB_{2p-2} = \lan \CO_\TZ,\CK^*_{|\TZ},S^2\CK^*_{|\TZ},\dots,S^{p-2}\CK^*_{|\TZ} \ran_{\D^b(Z)}, &
\text{for $n = 2p+1$.}
\end{array}
$$
Then $\D^b(\TZ) = \lan \CB_{n-3}(3-n),\CB_{n-4}(4-n),\dots,\CB_1(-1),\CB_0\ran$ is a $Z$-linear Lefschetz decomposition.
\end{lemma}
\begin{proof}
The restriction of the above collection to every fiber of $\TZ$ over $Z$ gives a Lefschetz
decomposition by~\cite{K3}. We conclude by~\cite{S}.
\end{proof}

We conclude that by theorem~\ref{th1} the category
$$
\TD = {}^\perp\Big\langle i_*(\CB_{n-3}(3-n)),\dots,i_*(\CB_1(-1))\Big\rangle \subset \D^b(\TY)
$$
is a categorical resolution of singularities of $Y$.

Let $H_G$ and $H_Y$ denote the divisor classes of a hyperplane sections
of $\Gr(n-4,W) \subset \PP(\Lambda^{n-4}W)$ and $Y \subset \PP(\Lambda^2W^*)$ respectively.

\begin{lemma}
We have $\TZ \sim 2H_Y - H_G$  and $K_\TY \sim (n-3)\TZ - 2nH_Y$.
In particular, $K_Y = -2nH_Y$.
\end{lemma}
\begin{proof}
Since $\TY = \PP_{\Gr(n - 4,W)}(\Lambda^2\CK^\perp)$, the Picard group
of $\TY$ is generated by $H_G$ and $H_Y$, hence $\TZ \sim \lambda H_G + \mu H_Y$
for some $\lambda,\mu \in \ZZ$. Moreover, it easy to see that $K_\TY = - (n-3)H_G - 6H_Y$.
Similarly, since $\TZ = \Fl(n-4,n-2;W)$, the Picard group
of $\TZ$ is generated by $H_G$ and $H_Y$, and $K_\TZ = - (n-2)H_G - 4H_Y$.
By adjunction formula we find $\lambda = -1$, $\mu = 2$.
Thus $K_\TY$ is linearly equivalent to $(n-3)\TZ - 2nH_Y$
which means that $K_Y = -2nH_Y$.
\end{proof}

We see that the discrepancy of $\TZ$ in $\TY$ equals to $n-3$, which is exactly by 1 less then
the number of the components of the Lefschetz decomposition of $\D^b(\TZ)$. Moreover, when $n$
is odd the Lefschetz decomposition is rectangular, hence the resolution $\TD$ of $Y$ is crepant
if $n$ is odd. Moreover, applying proposition~\ref{serre} we can find the Serre functor of $\TD$.

\begin{corollary}
The category $\TD$ has a Serre functor $\FS_\TD$. Moreover, if $F \in \TD$ is such that $i^*F \in \CB_{m-1}$,
then we have
$$
\FS_\TD(F) \cong F\otimes \pi^*\CO_Y(-2nH_Y)[4n-11].
$$
In particular, if $n$ is odd, the Serre functor is isomorphic to the tensoring by $\pi^*\CO_Y(-2nH_Y)[4n-11]$.
\end{corollary}

Finally, let us show that all these resolutions are noncommutative.
Let $E = \CO_\TY \oplus \CK^* \oplus \dots \oplus S^{\lfloor n/2\rfloor -2}\CK^*$, so that~\eqref{ldtz1} evidently holds.

\begin{lemma}
Vector bundle $E = \CO_\TY \oplus \CK^* \oplus \dots \oplus S^{\lfloor n/2\rfloor -2}\CK^*$ on $\TY$ is tilting over $Y$.
\end{lemma}
\begin{proof}
We use the same argument as in the proof of proposition~\ref{criterion}
replacing the projection to the exceptional divisor (which don't exists in our case)
by the projection $q:\TY = \PP_{\Gr(n-4,W)}(\Lambda^2\CK^\perp) \to \Gr(n-4,W)$.
We have to check that $R^{>0}\pi_*(\CEnd E) = 0$.
Since $Y$ is projective, it is equivalent to the equality
$H^{>0}(Y,\pi_*(\CEnd E)\otimes\CO_Y(t)) = 0$ for $t\gg0$.
By the projection formula we can rewrite this as
$H^{>0}(\TY,\CEnd E\otimes\pi^*\CO_Y(t))$.
But $\TY = \PP_{\Gr(n-4,W)}(\Lambda^2\CK^\perp)$.
Note that $E$ is a pullback from $\Gr(n-4,W)$.
Therefore, $q_*(\CEnd E\otimes\pi^*\CO_Y(t)) \cong \CEnd(\CO_{\Gr(n-4,W)}\oplus\CK^* \oplus \dots \oplus S^{\lfloor n/2\rfloor -2}\CK^*)\otimes S^t(\Lambda^2(W/\CK))$,
so it suffices to check that
$$
H^{>0}(\Gr(n-4,W),\CEnd(\CO_{\Gr(n-4,W)}\oplus\CK^* \oplus \dots \oplus S^{\lfloor n/2\rfloor -2}\CK^*)\otimes S^t(\Lambda^2(W/\CK))) = 0
\qquad\text{for $t\gg0$}.
$$
But the Borel--Bott--Weil theorem~\cite{D} easily implies vanishing of the above cohomology
for all $t\ge 0$.
\end{proof}

Thus, by theorem~\ref{th2}
$$
\TD = \D^b(Y,\CA)
$$
is a noncommutative resolution with $\CA = \pi_*\CEnd E$. If $n$ is odd this resolution is crepant.

\end{document}